\documentclass{amsart}

\usepackage{amssymb}
\usepackage{amsmath}
\usepackage{hyperref}
\usepackage{bm}

\newcommand{\half}{\frac{1}{2}}
\newcommand{\ep}{\epsilon}
\newcommand{\thalf}{\tfrac{1}{2}}
\newcommand{\summ}{\mathop{{\sum}^{P}}}

\newcommand{\sumh}{\mathop{{\sum}^{P}}}

\numberwithin{equation}{section}

\newtheorem{theorem}{Theorem}[section]

\newtheorem{lemma}[theorem]{Lemma}

\makeatletter
\@namedef{subjclassname@2010}{%
  \textup{2010} Mathematics Subject Classification}
\makeatother

\begin{document}

\title{The fifth moment of Hecke $L$-functions in the weight aspect}

\author{Rizwanur Khan}
\address{
Science Program\\ Texas A\&M University at Qatar\\ PO Box 23874\\ Doha, Qatar}
\email{rizwanur.khan@qatar.tamu.edu }


\begin{abstract}  
We prove an upper bound for the fifth moment of Hecke $L$-functions associated to holomorphic Hecke cusp forms of full level and weight $k$ in a dyadic interval $K\le k \le 2K$, as $K\to\infty$. The bound is sharp on Selberg's eigenvalue conjecture.
 \end{abstract}

\maketitle

\section{Introduction}

Moments of $L$-functions, especially at the central point, are extensively studied. They yield valuable data about an $L$-function's distribution, and can be used for example to infer information about the size, non-vanishing and symmetry-type of the central values.

This article is inspired by the recent works of Kiral and Young \cite{kiryou} and Blomer and Khan \cite{blokha}. The former paper established, for the first time, an upper bound for the fifth moment of $L$-functions associated to holomorphic newforms of prime level $q$ and fixed small weight, as $q\to \infty$. The latter paper established a certain reciprocity-type formula for the twisted fourth moment of Hecke $L$-functions in the level aspect, which gave as a corollary an upper bound for the fifth moment, but with more general conditions and also allowing for Maass $L$-functions. In both papers, the upper bound for the fifth moment depends on the Ramanujan conjecture at the finite places, and when assuming the truth of this conjecture, the given upper bound is sharp (as strong as the Lindel\"{o}f bound on average).

The goal of the present paper is to fix the level (at 1) and prove a fifth moment estimate in the weight aspect (it should also be possible to work with Hecke Maass $L$-functions in the spectral aspect). Let $H_k$ denote the orthonormal set of holomorphic Hecke cusp forms $f$ of level $1$
and weight $k$. This has $k/12+O(1)$ elements and forms a basis of the space of cusp forms of level 1 and weight $k$. Let $\lambda_f(n)$ denote the (real) eigenvalue corresponding to $f\in H_k$ of the $n$-th Hecke operator (which satisfies Deligne's bound $\lambda_f(n)\ll n^\ep$). The $L$-function associated to $f$ is defined for $\Re(s)>1$ by
\begin{align*}
L(s, f)=  \sum_{n=1}^{\infty}
\frac{\lambda_f(n)}{n^s}.
\end{align*}
The central point is $s=\half$ and by \cite{kohzag} the central value $L(\half,   f)$ is known to be non-negative. Our main theorem is 
\begin{theorem} \label{main} Let
\begin{align*}
\mathcal{F} = \bigcup_{\substack{K \le k \le 2K\\k\equiv 0 \bmod 2}} H_k,
\end{align*}
a set of $O(K^2)$ elements.
For any $\ep>0$, we have
\begin{align}
\label{fifthbound} \sum_{f\in \mathcal{F}} L(\thalf,f)^5 \ll K^{2+2\theta+\ep}
\end{align}
as $K\to \infty$, where $\theta= \frac{7}{64}$ is the current best bound towards the Selberg eigenvalue conjecture \cite[Appendix 2]{kim}.
\end{theorem}
\noindent The  ``log of conductor to log of family size'' ratio in (\ref{fifthbound}) is 5/2, the same as in the level aspect fifth moment considered in \cite{kiryou} and \cite{blokha}. Thus our result should be considered an analogue of the level aspect estimate. Assuming the Selberg eigenvalue conjecture (which is a part of the Ramanujan conjecture at the infinite place), our bound is sharp. This seems to be the first time that a sharp bound has been proven (conditionally) for any moment higher than the fourth in the archimedian (weight or spectral) aspect. Jutila \cite{jut} proved a good upper bound for the twelfth moment of Hecke Maass $L$-functions in the spectral aspect, but that is not sharp. 

Other authors \cite{ivi, jut2, pen} have proven sharp bounds for the third and fourth moments over smaller families. For example, in \cite{pen} Peng proved a sharp bound for the third moment over $H_k$, which yields the Weyl-quality bound $L(\half, f)\ll  k^{\frac13+\ep}$. Since such a strong bound already exists, we do not pursue a twisted fourth moment and amplification, although our methods would permit it. The goal is not to obtain individual bounds, although our main theorem already implies a weaker subconvexity bound.

Our main ideas have a similar flavour to those of \cite{kiryou,blokha,li}, but our method is different -- for example, we apply ``reciprocity'' twice, while the other papers apply it once. Compared to \cite{kiryou}, our proof is simpler and shorter, and as already noted above, our method could also be used to prove a bound for the twisted fourth moment, while this is not the case in \cite{kiryou} (because as explained in section 2 of that paper, the assumption $m_1\le m_2$ is made at the outset and cancellation in the $m_1$ sum is used to deal with ``fake main terms''). We cannot really compare with \cite{blokha} because that paper was after a more general result. It might be possible to derive our result from \cite{blokha} by first understanding the relevant integral transforms in terms of the weight, but our paper is self-contained and has the advantage (depending on taste) of being more ``classical'' in its approach. 

Throughout the paper, we will use the convention that $\ep$ denotes an arbitrarily small positive constant, but not necessarily the same one from one occurrence to the next.

\section{Rough Sketch}

The purpose of this sketch is to explain the main ideas, ignoring all technicalities. We will consider only the generic ranges of all sums.

Using approximate function equations, we can write the fifth moment as
\begin{align*}
\frac{1}{K^2} \sum_{f\in \mathcal{F}} L(\thalf,f)^5 &\approx \frac{1}{K^2} \sum_{f\in \mathcal{F}} \sum_{n_1\asymp K} \frac{\lambda_f(n_1)}{\sqrt{n_1}} \sum_{n_2,n_3,_4,n_5\asymp K} \frac{\lambda_f(n_2n_3n_4n_5)}{\sqrt{n_2n_3n_4n_5} }\\
&\approx \frac{1}{K^{\frac72}} \sum_{K\le k\le K} \frac{1}{K} \sum_{f\in H_k} \lambda_f(n_1) \lambda_f(n_2n_3n_4n_5).
\end{align*}
We need an upper bound of $O(K^{2\theta+\ep})$. We will in fact find that this kind of grouping with $n_1$ on one side and $n_2,n_3,n_4,n_5$ on the other leads to cleaner calculations. Applying Petersson's trace formula, the off-diagonal part of this is
\begin{align*}
 \frac{1}{K^{\frac72}} \sum_{K\le k\le K} \ \sum_{n_1,n_2,n_3,n_4,n_5 \asymp K} \  \sum_{c\ge 1} 2\pi i^k \frac{S(n_1, n_2n_3n_4n_5,c)}{c} J_{k-1} \Big(4\pi \frac{\sqrt{n_1n_2n_3n_4n_5}}{c}\Big).
\end{align*}
Summing over $k$ first, we will get that this is 
\begin{align*}
& \frac{1}{K^{\frac72}} \sum_{n_1,n_2,n_3,n_4,n_5 \asymp K} \sum_{c\asymp K^\half }  \frac{S(n_1, n_2n_3n_4n_5,c)}{c} e\Big( \frac{2\sqrt{n_1n_2n_3n_4n_5}}{c}\Big)\\
  \approx  &\frac{1}{K^{4}} \sum_{n_1,n_2,n_3,n_4,n_5 \asymp K} \sum_{c\asymp K^\half }  S(n_1, n_2n_3n_4n_5,c) e\Big( \frac{2\sqrt{n_1n_2n_3n_4n_5}}{c}\Big),
\end{align*}
where as usual $e(z)$ denotes $e^{2\pi i z}$.
Splitting the sum over $n_1$ into residue classes mod $c$ and applying Poisson summation (denote the dual variable by $m_1$) we get
\begin{align*}
  \frac{1}{K^{\frac72}}  \sum_{\substack{n_2,n_3,n_4,n_5 \asymp K\\ c\asymp K^\half} } \sum_{-\infty < m_1 <\infty} \sum_{a\bmod c} S(a, n_2n_3n_4n_5,c) e\Big( \frac{a m_1}{c}\Big) \int_{x\asymp 1} e\Big( \frac{2\sqrt{ x K n_2n_3n_4n_5}}{c}\Big) e\Big(\frac{-xK m_1}{c}\Big) dx.
\end{align*}
The complete sum over residue classes evaluates to $ce(\frac{-n_2n_3n_4n_5\overline{m_1}}{c})$, and the integral is evaluated using the stationary phase method. We get
\begin{align}
\label{return} \frac{1}{K^{4}}  \sum_{\substack{n_2,n_3,n_4,n_5\asymp K\\ c\asymp K^\half  \\  m_1\asymp K^\frac32 }}  e\Big( \frac{-n_2n_3n_4n_5 \overline{m_1}}{c}\Big) e\Big( \frac{n_2n_3n_4n_5}{m_1c}\Big) =  \frac{1}{K^{4}}  \sum_{\substack{n_2,n_3,n_4 ,n_5\asymp K\\ c\asymp K^\half  \\  m_1\asymp K^\frac32 }}  e\Big( \frac{n_2n_3n_4n_5\overline{c}}{m_1}\Big),
\end{align}
by reciprocity. 

Next we apply Poisson summation (mod $m_1$) to the $n_2$ and $n_3$ sums (in the actual proof, we will apply Voronoi summation once instead of Poisson summation twice).  Note that if we were following \cite{kiryou} step by step, we would have applied Poisson summation to $n_2$, $n_3$ and $n_4$, but this is not how we proceed. We get 
\begin{align}
\label{return2} \frac{1}{K^{\frac72}}  \sum_{\substack{n_4,n_5 \asymp K\\ c, m_2 ,m_3 \asymp K^\half  \\  m_1\asymp K^\frac32 }}  e\Big( \frac{m_2m_3 c \overline{n_4n_5}}{m_1}\Big).
\end{align}
This sum displays only the generic ranges of $m_2$ and $m_3$ (the dual variables). The zero frequencies $m_2=0$ or $m_3=0$, which are omitted, are in fact quite troublesome. For example, return to (\ref{return}) and consider the terms with $m_1|n_3n_4n_5$ (these terms correspond to $m_2=0$). The contribution of such terms is 
\begin{align}
\label{fakelarge} \frac{1}{K^{4}}  \sum_{\substack{n_2,n_3,n_4,n_5 \asymp K\\ c\asymp K^\half , m_1\asymp K^\frac32 \\ m_1| n_3n_4n_5}}  1 \asymp K^{\half},
\end{align}
while we need to prove a bound of $K^{2\theta+\ep}$. It seems that we cannot do better because there are no harmonics present to produce further cancellation. Of course, it is not possible (by the Lindel\"{o}f hypothesis) for the fifth moment to be so large, so a careful evaluation of the fifth moment must show that these ``fake main terms'' should cancel out somehow. But there is a shortcut. The weight functions from the approximate functional equations have been suppressed in (\ref{fakelarge}). If we take them into account, there is a way to design them carefully so that (\ref{fakelarge}) is not so large. This idea was used in \cite{bhm} and \cite{kiryou}, and section 2 of the latter paper contains a nice heuristic about how the idea works.  

Back to (\ref{return2}), we can apply reciprocity again to get
\begin{align*}
 \frac{1}{K^{\frac72}}  \sum_{\substack{n_4,n_5 \asymp K\\ c, m_2 ,m_3 \asymp K^\half  \\  m_1\asymp K^\frac32 }}  e\Big( \frac{m_2m_3 c \overline{m_1}}{n_4n_5}\Big)  e\Big( \frac{-m_2m_3 c}{m_1n_4n_5}\Big) \approx  \frac{1}{K^{\frac72}}  \sum_{\substack{n_4,n_5 \asymp K\\ c, m_2 ,m_3 \asymp K^\half  \\  m_1\asymp K^\frac32 }}  e\Big( \frac{m_2m_3 c \overline{m_1}}{n_4n_5}\Big).
\end{align*}
Applying Poisson summation (mod $n_4n_5$) to the $m_1$ sum (denote the dual variable by $l_1$), this is
\begin{align*}
 \frac{1}{K^{4}}  \sum_{\substack{n_4,n_5 \asymp K\\ c, m_2 ,m_3 , l_1 \asymp K^\half }} S(m_2m_3 c , l_1 ,n_4n_5).
\end{align*}
Now we can sum over $n_4$ using Kuznetsov's formula. The sum of Kloosterman sums is in the Linnik range as $n_4n_5\ge \sqrt{ m_2 m_3 c l_1}$. This leads to 
\begin{align}
\label{sketchlast} \frac{1}{K^{2}}  \sum_{n_5 \asymp K}  \Big( \sum_{\substack{ c, m_2 ,m_3 , l_1 \asymp K^\half }} \sum_{t_j\asymp 1 } \frac{\lambda_j(m_2m_3 c) \lambda_j(l_1)}{\sqrt{m_2 m_3 c l_1}} + \ldots \Big)
\end{align}
where the sum is over an orthonormal basis of Maass cusp forms $\{ u_j \}$ of level $n_5$ and (essentially bounded) Laplacian eigenvalue $\frac{1}{4}+t_j^2$, and the ellipsis denotes the contribution of the Eisenstein series and holomorphic forms. Actually we lose $O(K^{2\theta+\ep})$ here due to the possibility of exceptional eigenvalues, but for the purposes of this sketch we ignore this issue. 

The inner sum of (\ref{sketchlast}), given within the parentheses, looks like the fourth moment of $L(\half, u_j)$ in the level aspect, provided that we can decompose $\lambda_j(m_2m_3 c)$ by multiplicativity. For this, we need to work with a basis comprising of lifts of newforms; such a basis is given in \cite{blomil} or \cite{blokha}. Then the expected bound for the fourth moment, which can be proved using the spectral large sieve, gives 
\begin{align*}
 \frac{1}{K^{2}}  \sum_{n_5 \asymp K} (n_5 K^\ep) \ll K^\ep
\end{align*}
as desired. We never need any cancellation from the $n_5$-sum, which is why a twisted fourth moment bound would probably be possible in place of the main theorem.

\section{Background}

\subsection{Approximate functional equations}

For $f\in H_k$ we have the functional equation \cite[Theorem 14.7]{iwakow},
\begin{align}
\label{fe} \Lambda(s,f):= (2\pi)^{-s} \Gamma(s+\tfrac{k-1}{2})L(s,f) = i^k \Lambda(1-2,f).
\end{align} 
Let $\tau(m)$ denote the number of divisors of $m$. We will use the following standard approximate functional equations. For any $f\in H_k$, we have

\begin{align}
\label{afe1} L(\thalf,f)^2 = 2\sum_{m\ge 1} \frac{\lambda_f(m)\tau(m)}{\sqrt{m}} V_k(m),
\end{align}
where
\begin{align*}
V_k(x)=\frac{1}{2\pi i} \int_{(A)} x^{-s} \mathcal{G}(s)  \frac{\Gamma(s+\tfrac{k}{2})^2}{\Gamma(\tfrac{k}{2})^2} \zeta(1+2s) \frac{ds}{s}
\end{align*}
for any $A>0$ and
\begin{align}
\label{gdef} \mathcal{G}(s)= 2 e^{s^2} (\thalf-s^2).
\end{align}
This follows from the functional equation (\ref{fe}) and \cite[Theorem 5.3]{iwakow}. As explained in that theorem, we may insert in the integrand above any even function which is bounded in a fixed horizontal strip about $\Re(s)=0$, and has value 1 at $s=0$. Our function $\mathcal{G}(s)$ satisfies these properties and is chosen to decay exponentially in the vertical direction (this is convenient for convergence) and to vanish at $s=\half$ (this will be needed later to deal with the ``fake main terms'').

For $k\equiv 0 \bmod 4$, the root number in the functional equation is $1$, and we have
\begin{align}
\label{afe2} L(\thalf,f) = 2\sum_{n\ge 1} \frac{\lambda_f(n)}{\sqrt{n}} W_k(n),
\end{align}
where
\begin{align*}
W_k(x)=\frac{1}{2\pi i} \int_{(A)} x^{-s}   e^{s^2} \frac{\Gamma(s+\tfrac{k}{2})}{\Gamma(\tfrac{k}{2})} \frac{ds}{s}.
\end{align*}

We have
\begin{align}
\label{vwbounds} V_k^{(j)}(x), W_k^{(j)}(x) \ll x^{-j}(1+x)^{-A}
\end{align}
for any $A>0$ and integer $j\ge 0$.
Using this for $j=0$, large $A$ and Stirling's estimates for the gamma function, the sums (\ref{afe1}) and (\ref{afe2}) may be restricted to $m\ll k^{2+\ep}$ and $n\ll k^{1+\ep}$ respectively, up to an error of $O(k^{-100})$.  Taking $j=0$ and $A=\ep$ shows that $|V_k(x)|, |W_k(x)| <k^\ep$.

\subsection{Summation formulae}

 We will need the Voronoi summation formula and the Poisson summation formula.
 
 \begin{lemma} \label{vorsum} {\bf Voronoi summation.}
Given a compactly supported smooth function $\Phi$ with bounded derivatives, and coprime integers $h$ and $\ell$, we have
\begin{align}
\label{vorstate} \sum_{m\ge 1} \frac{\tau(m)}{m} e\Big(\frac{m\overline{h}}{\ell}\Big)\Phi\Big(\frac{m}{M}\Big)= \frac{1}{\ell} \int_{-\infty}^\infty \log\Big(\frac{ x}{\ell^2}+2\gamma\Big)\Phi \Big(\frac{x}{M}\Big)\frac{dx}{x} + \sum_\pm \frac{1}{\ell} \sum_{r\ge 1} \tau(r) e\Big(\frac{\pm rh}{\ell}\Big) \check{\Phi}_\pm \Big(\frac{Mr}{\ell^2}\Big),
\end{align}
where
\begin{align*}
\check{\Phi}_\pm(x)= \frac{1}{2\pi i} \int_{(A)} H_1^\pm(s) \tilde{\Phi}(-s)  x^{-s}  ds,
\end{align*}
$\tilde{\Phi}$ is the Mellin transform of $\Phi$, 
\begin{align*}
H_1^\pm(s)= 2(2\pi i)^{-2s} \Gamma(s)^2  \cos^{(1\mp1)/2} (\pi s),
\end{align*}
and $A>0$. 
\end{lemma}
\proof 
See \cite[section 2.3]{bhm}. We can take any $A>0$ because $\tilde{\Phi}(-s)\ll (1+|s|)^{-B}$ for any $B\ge 0$ by integration by parts. 
\endproof

\begin{lemma} \label{poiss} {\bf Poisson summation.} Given a compactly supported smooth function $\Phi$ with bounded derivatives, and an arithmetic function $S_q(n)$ with period $q$, we have
\begin{align}
\label{pois3} &\sum_{-\infty < n < \infty } \Phi\Big(\frac{n}{N}\Big) S_q(n) \\
\nonumber &= \frac{N}{q} \sum_{-\infty < l <\infty} \hat{\Phi}\Big(\frac{lN}{q}\Big) \sum_{a\bmod q} S_q(a) e\Big(\frac{a l}{q}\Big)\\
\nonumber &=  \frac{N}{q}  \hat{\Phi} (0)  \sum_{a\bmod q} S_q(a) + \frac{N}{q} \sum_{-\infty < l <\infty} \  \sum_{a\bmod q} S_q(a) e\Big(\frac{a l}{q}\Big) \int_{-\infty}^{\infty} \frac{1}{2\pi i} \int_{(A)} \Big(\frac{-2\pi x lN}{q}\Big)^{-s} \Phi(x) H_2(s)  ds dx,
\end{align}
where $\hat{\Phi}$ denotes the Fourier transform of $\Phi$,
\begin{align*}
H_2(s)=   \Gamma(s) \exp \Big(\frac{i\pi s}{2}\Big) 
\end{align*}
and $A>0$.
\end{lemma}
\proof
For the second line of (\ref{pois3}), separate the $n$ sum into sums over residue classes $a$ modulo $q$ and apply the usual Poisson summation formula to each sum. For the third line we keep aside the contribution of $l=0$, and for $l\neq 0$ we first compute the Mellin transform
\begin{align}
\label{mel} \int_0^\infty \hat{\Phi}(y) y^{s-1} dy = \int_0^\infty \int_{-\infty}^\infty \Phi(x) e(-yx) y^{s-1} dx dy =\int_{-\infty}^\infty \Phi(x) (-2\pi x)^{-s} H_2(s) dx.
\end{align}
This follows by swapping the order of integration, which we can do by the compact support of $\Phi$, and then using the Mellin transform 
\begin{align*}
\int_0^\infty e^{iy} y^{s-1} dy = H_2(s)
\end{align*}
which holds for $0<\Re(s)<1$. But since $\hat{\Phi}(y)\ll (1+|y|)^{-B}$ for any $B\ge 0$ by integration by parts, we have that $\int_0^\infty \hat{\Phi}(y) y^{s-1} dy $ converges absolutely for $\Re(s)>0$. Thus the Mellin transform given in (\ref{mel}) analytically continues to $\Re(s)>0$, and by the Mellin inversion theorem we have
\begin{align*}
\hat{\Phi}\Big(\frac{lN}{q}\Big)= \frac{1}{2\pi i} \int_{(A)} \Big(\frac{lN}{q}\Big)^{-s} \int_{-\infty}^\infty \Phi(x) (-2\pi x )^{-s} H_2(s) dx ds
\end{align*}
for any $A>0$.
\endproof

\subsection{An average of the $J$-Bessel function}  The following result can be found in \cite[Corollary 8.2]{ils}.
\begin{lemma}  \label{javg}
Let $x>0$ and let $h$ be a smooth function compactly supported on the positive reals and possessing bounded derivatives. We have
\begin{align}
\label{avgmaint} \frac{1}{K} \sum_{k \equiv 0 \bmod 2} 2i^k
h\Big(\frac{k-1}{K}\Big)
J_{k-1}(x)
 = -\frac{1}{\sqrt{x}} \Im \Big(e^{-2\pi i /8} e^{ix}
\hbar\Big(\frac{K^2}{2x}\Big) \Big) +
O\Big(\frac{x}{K^5}\int_{-\infty}^{\infty} v^4|\hat{h}(v)|dv\Big),
\end{align}
where for real $v$,
\begin{align*}
\hbar(v):=\int_{0}^{\infty} \frac{h(\sqrt{u})}{\sqrt{2\pi u}}
e^{iuv} du
\end{align*}
 and $\hat{h}$ denotes the Fourier transform of $h$.
The implied constant is absolute.
\end{lemma}
\noindent By integrating by parts
several times we get that $\hbar(v)\ll |v|^{-B}$
 for any $B\ge 0$. Thus the main
term of (\ref{avgmaint}) is not dominant if $x < K^{2-\ep}$.

For future use, define for any complex number $s$ the more general function
\begin{align*}
\hbar_{s}(v):=\int_{0}^{\infty}
\frac{h(\sqrt{u})}{\sqrt{2\pi u}} u^{s/2} e^{iuv} du.
\end{align*}
Integrating by parts, we get 
\begin{align}
\label{hbar} \hbar_{s}^{(j)}(v)\ll_{\Re(s)} (1+|s|)^B |v|^{-B}
\end{align}
 for any $B\ge 0$. Thus the Mellin transform
 \begin{align*}
\tilde{\hbar}_{s}(w)= \int_0^\infty \hbar_{s}(v) v^{w-1} dv 
\end{align*}
is holomorphic in the half plane $\Re(w)>0$, and we have by integrating by parts $j$ times:
\begin{align*}
\tilde{\hbar}_{s}(w) \ll_{\Re(s)} (1+|s|)^{j+\Re(w)+1} (1+|w|)^{-j}.
\end{align*}

\section{Hecke relations}

Define
\begin{align*}
\summ_{f\in H_k} \gamma_f := \sum_{f\in H_k}
\Big(\frac{2\pi^2}{(k-1)L(\mathrm{1, sym}^2 f )}\Big) \gamma_f
\end{align*}
for any complex numbers $\gamma_f$ depending on $f$. The average $\summ$ arises in the Petersson trace formula \cite[Proposition 14.5]{iwakow}:
\begin{align*}
\summ_{f\in H_k} \lambda_{f}(n) \lambda_f(m) =
\delta_{m,n}+2\pi i^k \sum_{c=1}^{\infty}
\frac{S(n,m,c)}{c}J_{k-1}\Big(\frac{4\pi\sqrt{mn}}{c}\Big),
\end{align*}
where the value of $\delta_{m,n}$ is $1$ if $m=n$ and $0$ otherwise,
$S(n,m,c)$ is the Kloosterman sum, and $J_{k-1}(x)$ is the $J$-Bessel
function.

The following lemma explains how we will group together variables in the fifth moment.

\begin{lemma} \label{reduce}
To prove the main theorem, it suffices to prove that for any smooth functions $h, U_1,U_2,U_3$ compactly supported on $(\half,\frac{5}{2})$ with bounded derivatives, and any
\begin{align*}
\alpha,\beta,\beta_1,\beta_2\ge 1, \ \ 1\le N_1,N_2,N_3 < K^{1+\ep},  
\end{align*}
with
\begin{align}
\label{condition} N_3\ge N_2, \ \ N_1N_2<\frac{K^{2+\ep}}{\alpha}, \ \  \beta\ge \alpha,
\end{align}
we have
\begin{align}
\label{hrel} 
\frac{1}{K} \sum_{k \equiv 0 \bmod 2} h\Big(\frac{k-1}{K}\Big) \sumh_{f \in H_k} S_f  \ll \sqrt{\alpha} K^{2\theta+\ep},
\end{align}
where
\begin{multline*}
S_f :=\\\sum_{n_1,n_2,n_3,m\ge 1} \frac{\lambda_f(n_1 n_2 m \alpha ) \lambda_f(n_3) \tau(m) }{ \sqrt{n_1 n_2 n_3m }} 
W_k(n_1\beta_1)W_k(n_2\beta_2)W_k(n_3) V_k(m\beta) U_1\Big(\frac{n_1 }{N_1}\Big)U_2\Big(\frac{ n_2 }{N_2}\Big)U_3\Big(\frac{n_3}{N_3}\Big).
\end{multline*}
\end{lemma}
\proof
To prove the main theorem, it suffices to prove that
\begin{align*}
\frac{1}{K} \sum_{k \equiv 0 \bmod 2} h\Big(\frac{k-1}{K}\Big) \sumh_{f \in H_k} L(\thalf,f)^5  \ll K^{2\theta+\ep},
\end{align*}
because we have $L(\half,f)\ge 0$ by \cite{kohzag} and $k^{-\ep}<L(\mathrm{1, sym}^2 f )<k^{\ep}$ by \cite[Appendix]{golhoflie}.

We claim that
\begin{align}
\label{claim} L(\thalf,f)^5 &= 8 \Big( \sum_{n \ge 1} \frac{\lambda_f(n)}{\sqrt{n}} W_k(n) \Big)^3 L(\thalf,f)^2.
\end{align}
This holds by (\ref{afe2}) when $k\equiv 0 \bmod 4$. But when $k\equiv 2 \bmod 4$, it also holds because then $L(\half,f)=0$ by the functional equation (\ref{fe}), so both sides of (\ref{claim}) vanish. Now we can insert the approximate functional equation for $L(\thalf,f)^2$ given in (\ref{afe1}) to get that
\begin{align*}
L(\thalf,f)^5 &= 16 \Big( \sum_{n \ge 1} \frac{\lambda_f(n)}{\sqrt{n}} W_k(n) \Big)^3 \Big(\sum_{m \ge 1} \frac{\lambda_f(m )\tau(m )}{\sqrt{m}} V_k(m) \Big).
\end{align*}
Expanding the cube and working in dyadic intervals, to establish the main theorem it suffices to prove that 
\begin{align*}
\frac{1}{K} \sum_{k \equiv 0 \bmod 2} h\Big(\frac{k-1}{K}\Big) \sumh_{f \in H_k} S_1  \ll K^{2\theta+\ep},
\end{align*}
where
\begin{align*}
S_1:= \prod_{i=1}^3\Big( \sum_{n_i \ge 1} \frac{\lambda_f(n_i)}{\sqrt{n_i}} W_i\Big(\frac{n_i}{N_i}\Big) \Big) \Big(\sum_{m \ge 1} \frac{\lambda_f(m )\tau(m )}{\sqrt{m}} V_k(m) \Big)
\end{align*}
for
\begin{align*}
W_i(n_i):=W_k(n_i)U_i\Big(\frac{n_i}{N_i}\Big)
\end{align*}
and $1\le N_1,N_2,N_3<K^{1+\ep}$. By symmetry, we can suppose that $N_3\ge N_2$. By Hecke multiplicativity, we have
\begin{align*}
 \lambda_f(m)\lambda_f(n_1)=\sum_{d|(m,n_1)} \lambda_f\Big(\frac{mn_1}{d^2}\Big),
\end{align*}
so replacing $m$ by $md$ and $n_1$ by $n_1d$, we get
\begin{align*}
S_1= \sum_{n_1,n_2,n_3,m,d\ge 1} \frac{\lambda_f(mn_1)\lambda_f(n_2)\lambda_f(n_3)\tau(md)}{d\sqrt{n_1n_2n_3m}} W_1\Big(\frac{n_1d}{N_1}\Big)W_2\Big(\frac{n_2}{N_2}\Big)W_3\Big(\frac{n_3}{N_3}\Big)V_k(md).
\end{align*}
Now we combine
\begin{align*}
\lambda_f(mn_1) \lambda_f(n_2) =\sum_{b|(mn_1,n_2)} \lambda_f\Big(\frac{mn_1n_2}{b^2}\Big)=\sum_{\substack{ n_2= b_1b \\ b|mn_1}} \lambda_f\Big(\frac{mn_1b_1}{b}\Big).
\end{align*}
Ordering by the gcd of $n_1$ and $b$, we have the disjoint union
\begin{align}
\label{proc1} \{ n_1,m :  b |n_1m \}  = \bigsqcup_{\substack{ b=b_2 b'  \\ (b,n_1)=b_2 }}  \{ n_1,m :  b |n_1m\}  =   \bigsqcup_{\substack{ b=b_2b'  }}   \Big\{ n_1,m :  b_2 | n_1, b' |m, \Big(\frac{n_1}{b_2},b'\Big)=1 \Big\},
\end{align}
and $(\frac{n_1}{b_2},b')=1$ can be detected using the Mobius function:
\begin{align}
\label{proc2} \sum_{\substack{b'=b_3b_4\\b_3|\frac{n_1}{b_2}}} \mu(b_3) = \begin{cases}
1 &\text{ if } (\frac{n_1}{b_2},b')=1\\
0 &\text{otherwise}.
\end{cases}
\end{align}
Thus replacing $b$ by $b_2b_3b_4$, $n_2$ by $b_1b_2b_3b_4$, $n_1$ by $n_1b_2b_3$, and $m$ by $mb_3b_4$, we get
\begin{multline*}
S_1= \sum_{\substack{n_1,b_1,n_3,m\ge 1\\ b_2,b_3,b_4,d\ge 1}} \frac{\lambda_f(mn_1 b_1 b_3) \lambda_f(n_3)\tau(mb_3b_4d) \mu(b_3) }{d b_2 b_3^\frac32 b_4 \sqrt{n_1b_1 n_3m}} \\W_1\Big(\frac{n_1 b_2 b_3d}{N_1}\Big)W_2\Big(\frac{b_1b_2b_3b_4}{N_2}\Big)W_3\Big(\frac{n_3}{N_3}\Big)V_k(mb_3b_4d).
\end{multline*}
Splitting the divisor function 
\begin{align*}
\tau(mb_3b_4d)=\sum_{r|(m,b_3b_4d)} \mu(r) \tau\Big(\frac{m}{r}\Big) \tau\Big(\frac{b_3b_4 d}{r}\Big),
\end{align*}
replacing $m$ by $mr$, and renaming $b_1$ to $n_2$, we have
\begin{multline*}
S_1= \sum_{\substack{n_1,n_2,n_3,m\ge 1 \\ b_2,b_3,b_4,d\ge 1\\ r| b_3b_4 d }} \frac{\lambda_f(mn_1 n_2 b_3 r) \lambda_f(n_3)\tau(m)\tau(\frac{b_3b_4 d}{r}) \mu(b_3) \mu(r) }{d b_2 b_3^\frac32 b_4 \sqrt{n_1 n_2 n_3m r}} \\
W_1\Big(\frac{n_1 b_2 b_3d}{N_1}\Big)W_2\Big(\frac{ n_2 b_2b_3b_4}{N_2}\Big)W_3\Big(\frac{n_3}{N_3}\Big)V_k(mb_3b_4dr).
\end{multline*}
We plan to find cancellation in the sum over $n_1,n_2,n_3,m$ and to sum trivially over the remaining parameters $b_2,b_3,b_4,r,d$. Thus it suffices to prove that
\begin{align}
\label{need0} \sum_{\substack{b_2,b_3,b_4,d\ge 1 \\ r| b_3b_4 d }}  \frac{1}{d b_2 b_3^\frac32 b_4 \sqrt{r}} \Bigg|\frac{1}{K} \sum_{k \equiv 0 \bmod 2} h\Big(\frac{k-1}{K}\Big) \sumh_{f \in H_k} S_2 \Bigg|  \ll K^{2\theta+\ep},
\end{align}
where
\begin{align*}
 S_2:=\sum_{n_1,n_2,n_3,m\ge 1} \frac{\lambda_f(mn_1 n_2 b_3 r) \lambda_f(n_3)\tau(m) }{ \sqrt{n_1 n_2 n_3m }} 
W_1\Big(\frac{n_1 b_2 b_3d}{N_1}\Big)W_2\Big(\frac{ n_2 b_2b_3b_4}{N_2}\Big)W_3\Big(\frac{n_3}{N_3}\Big)V_k(mb_3b_4dr).
\end{align*}
For (\ref{need0}) it suffices to show that
\begin{align*}
\frac{1}{K} \sum_{k \equiv 0 \bmod 2} h\Big(\frac{k-1}{K}\Big) \sumh_{f \in H_k} S_2\ll K^{2\theta+\ep} \sqrt{b_3 r}.
\end{align*}
This is given by (\ref{hrel}), once in $S_f$ we replace $N_1$ by $N_1/b_2b_3d$ and $N_2$ by $N_2/b_2b_3b_4$, and take $\alpha=b_3 r$, $\beta=b_3b_4 d r$, $\beta_1= b_2b_3d$, $\beta_2= b_2b_3b_4$. Note that these substitutions lead to a smaller value of $N_2$, so that $N_3\ge N_2$ still holds. Since $\beta_1\beta_2\ge b_3(b_3b_4r)\ge \alpha$, we have $N_1N_2<K^{2+\ep}/\alpha$.  Also note that $\beta\ge \alpha$.
\endproof

\section{Application of the trace formula}

Applying the Petersson trace formula to Lemma \ref{reduce}, we need to prove that
\begin{align*}
D +OD< \sqrt{\alpha}K^{2\theta+\ep},
\end{align*}
where the diagonal
\begin{multline*}
D:=\\ \sum_{\substack{n_1,n_2,n_3, m \ge 1\\ n_3=n_1n_2m\alpha}} \frac{\lambda_f(n_1 n_2 m \alpha ) \lambda_f(n_3) \tau(m) }{ \sqrt{n_1 n_2 n_3m }} 
U_1\Big(\frac{n_1 }{N_1}\Big)U_2\Big(\frac{ n_2 }{N_2}\Big)U_3\Big(\frac{n_3}{N_3}\Big)W_k(n_1\beta_1)W_k(n_2\beta_2)W_k(n_3)V_k(m\beta)
\end{multline*}
trivially satisfies the required bound, and the off-diagonal is
\begin{multline*}
OD:= \sum_{\substack{n_1,n_2,n_3, m,c \ge 1\\ }} \frac{S(n_1 n_2 m \alpha , n_3 ,c) \tau(m) }{c \sqrt{n_1 n_2 n_3m }} U_1\Big(\frac{n_1 }{N_1}\Big)U_2\Big(\frac{ n_2 }{N_2}\Big)U_3\Big(\frac{n_3}{N_3}\Big) \\
\frac{1}{K}  \sum_{k \equiv 0 \bmod 2} h\Big(\frac{k-1}{K}\Big) 2\pi i^k J_{k-1}\Big(\frac{4\pi  \sqrt{n_1n_2n_3 m \alpha}}{c} \Big) W_k(n_1\beta_1) W_k(n_2\beta_2) W_k(n_3)V_k(m\beta).
\end{multline*}
At this point, we cannot absorb the $W_k$ functions into the arbitrary weight functions $U_i$ because $W_k$ depends on $k$ and we still need to average over $k$, which is what we do next. Applying Lemma \ref{javg}, the contribution of its error term is bounded by 
\begin{align*}
\frac{1}{K^{5-\ep}} \sum_{\substack{n_1n_2<K^{2+\ep}/\alpha \\ n_2<K^{1+\ep}\\ m<K^{2+\ep}  }} \  \sum_{c\ge 1} \  \frac{|S(n_1 n_2 m \alpha , n_3 ,c)|  }{c \sqrt{n_1 n_2 n_3m }} \frac{ \sqrt{n_1 n_2 n_3m \alpha}}{c} \ll  K^\ep,
\end{align*}
on using Weil's bound for the Kloosterman sum. Thus we need only consider the main term of Lemma \ref{javg}, and it suffices to prove
\begin{align*}
OD_1 :=\sum_{\substack{n_1,n_2,n_3, m,c \ge 1\\ }}& \frac{S(n_1 n_2 m \alpha , n_3 ,c) \tau(m) }{\sqrt{c} (n_1 n_2 n_3m)^\frac34 } e\Big(\frac{2\pi  \sqrt{n_1n_2n_3 m \alpha}}{c} \Big) \\ &\Psi_K\Big(n_1\beta_1,n_2\beta_2,n_3, m\beta , \frac{K^2 c}{8\pi  \sqrt{n_1n_2n_3 m \alpha}}\Big)\prod_{j=1}^{3} U_j\Big(\frac{n_j}{N_j}\Big) \ll \alpha^\frac34 K^{2\theta+\ep},
\end{align*}
where
\begin{multline*}
\Psi_K(x_1,x_2,x_3, x_4, v) := \frac{1}{(2\pi i)^4} \int_0^\infty    \frac{h(\sqrt{u})}{\sqrt{2\pi u}} \int_{(A_4)}  \frac{ \zeta(1+2s_4) \mathcal{G}(s_4) }{ x_4^{s_4}}  \frac{\Gamma^2\Big(\frac{\sqrt{u}}{2}K+s_4+\frac{1}{2}\Big)}{\Gamma^2\Big(\frac{\sqrt{u}}{2}K+\frac{1}{2}\Big)} \\
\prod_{j=1}^{3} \int_{(A_j)}  \frac{e^{s_j^2} }{x_j^{s_j} } 
\frac{\Gamma\Big(\frac{\sqrt{u}}{2}K+s_j+\frac{1}{2}\Big)}{\Gamma\Big(\frac{\sqrt{u}}{2}K+\frac{1}{2}\Big)} \frac{ds_j}{s_j}\frac{ds_4}{s_4} e^{iuv}  du.
\end{multline*}
By the rapid decay of the $s_1,s_2,s_3,s_4$ integrands in vertical lines, we may effectively truncate the integrals to $|\Im s_1|,|\Im s_2|,|\Im s_3|,|\Im s_4|<K^\ep$.  For $|s|<K^\ep$, by Stirling's approximation we have
\begin{align*}
\frac{\Gamma\Big(\frac{\sqrt{u}}{2}K+s+\frac{1}{2}\Big)}{\Gamma\Big(\frac{\sqrt{u}}{2}K+\frac{1}{2}\Big)} = \Big(\frac{\sqrt{u}}{2} K \Big)^s\Big( 1 + \frac{P(s)}{\sqrt{u}K}+ O\Big(\frac{1}{K^{2-\ep}}\Big)\Big)
\end{align*}
for some polynomial $P$. Thus
\begin{align}
\label{wstir} \Psi_K(x_1,x_2,x_3,x_4,v)  = \Psi\Big(\frac{x_1}{K},\frac{x_2}{K},\frac{x_3}{K},\frac{x_4}{K^2},v\Big) + \frac{1}{K} \Psi_0\Big(\frac{x_1}{K},\frac{x_2}{K},\frac{x_3}{K},\frac{x_4}{K^2},v\Big) + O(K^{-2+\ep}),
\end{align}
where for $\xi_i>0$ and real $v$ we define
\begin{multline}
\label{wdef} \Psi(\xi_1,\xi_2,\xi_3,\xi_4,v) \\ := \frac{1}{(2\pi i)^4} \int_{(A_4)} \int_{(A_3)} \int_{(A_2)}  \int_{(A_1)}   \frac{e^{s_1^2+s_2^2+s_3^2} \zeta(1+2s_4) \mathcal{G}(s_4) }{ 2^{s_1+s_2+s_3} 4^{s_4}  \xi_1^{s_1} \xi_2^{s_2} \xi_3^{s_3} \xi_4^{s_4}} \hbar_{s_1+s_2+s_3+2s_4}(v)  \frac{ds_1}{s_1}\frac{ds_2}{s_2}\frac{ds_3}{s_3}\frac{ds_4}{s_4}
\end{multline}
and $W_0$ has the same definition except for the presence of an extra factor $\frac{P(s_1,s_2,s_3,s_4)}{\sqrt{u}}$ in the integrand for some polynomial $P$. It suffices to treat only the contribution of $\Psi$, as the treatment of the secondary term $\Psi_0$ will be similar. Thus we need to prove
\begin{align}
\label{defod2} OD_2 :=\sum_{\substack{n_1,n_2,n_3, m,c \ge 1\\ }}& \frac{S(n_1 n_2 m \alpha , n_3 ,c) \tau(m) }{\sqrt{c} (n_1 n_2 n_3m)^\frac34 } e\Big(\frac{2\pi  \sqrt{n_1n_2n_3 m \alpha}}{c} \Big) \\ 
\nonumber &\Psi\Big(\frac{n_1\beta_1}{K},\frac{n_2\beta_2}{K},\frac{n_3}{K}, \frac{m\beta }{K^2}, \frac{K^2 c}{8\pi  \sqrt{n_1n_2n_3 m \alpha}}\Big)\prod_{j=1}^{3} U_j\Big(\frac{n_j}{N_j}\Big) \ll \alpha^\frac34 K^{2\theta+\ep}.
\end{align}
By (\ref{hbar}) we may assume (up to negligible error) that
\begin{align}
\label{crange}  c\ll \frac{\sqrt{n_1 n_2 n_3 m \alpha}}{K^{2-\ep}}.
\end{align}
By (\ref{vwbounds}) and (\ref{hbar}), we have that
\begin{align}
\label{wbound} \frac{\partial^{j_1}}{\partial \xi_1^{j_1} }\frac{\partial^{j_2}}{\partial \xi_2^{j_2} }\frac{\partial^{j_3}}{\partial \xi_3^{j_3} } \frac{\partial^{j_4}}{\partial \xi_4^{j_4} } \frac{\partial^{j}}{\partial v^{j} } \Psi(  \xi_1,  \xi_2,  \xi_3,  \xi_4,   v) \ll K^\ep \xi_1^{-j_1-A_1} \xi_2^{-j_2-A_2} \xi_3^{-j_3-A_3} \xi_4^{-j_4-A_4} v^{-j-B}
\end{align}
for $\xi_1,\xi_2,\xi_3,\xi_4,v>0$, any integers $j_i,B\ge 0$ and any real $A_i>0$.

\section{Poisson summation and reciprocity}

In (\ref{defod2}), we sum over $n_3$ in residue classes mod $c$ and apply Poisson summation (Lemma \ref{poiss}), getting
\begin{multline}
\label{od2} OD_2 = \sum_{\substack{n_1,n_2,m,c \ge 1\\ -\infty < \ell < \infty} } \frac{N_3}{c} \sum_{a \bmod c} S(a,n_1 n_2 m \alpha,c) e\Big(\frac{a\ell}{c}\Big)  \frac{\tau(m)}{ \sqrt{c} (n_1 n_2 N_3m)^\frac34 } U_1\Big(\frac{n_1}{N_2}\Big)U_2\Big(\frac{n_2}{N_2}\Big)\\ \int_{-\infty}^{\infty} e\Big(\frac{-\ell N_3 x}{c}\Big) e\Big(\frac{2\sqrt{x n_1n_2N_3 m \alpha}}{c}\Big) \Psi\Big(\frac{n_1\beta_1}{K},\frac{n_2\beta_2}{K},\frac{xN_3}{K}, \frac{m\beta }{K^2}, \frac{K^2 c}{8\pi  \sqrt{x n_1n_2N_3 m \alpha}}\Big)  \frac{U_3(x)}{x^\frac34} dx.
\end{multline}
Call the integral above $I$. We will evaluate it using stationary phase approximation.
\begin{lemma} \label{stationary} We have that $I\ll K^{-1000}$ unless $|\ell|\le K^{100}$, in which case
\begin{align*}
I &=\\
  &\sqrt{\frac{2\ell c}{ n_1n_2 m \alpha}}  e\Big(\frac{n_1 n_2m\alpha }{\ell c}- \frac{\pi}{8}\Big) \Psi \Big(\frac{n_1\beta_1}{N_1},\frac{n_2\beta_2}{N_2},
\frac{ n_1n_2m\alpha}{\ell^2 K}, \frac{m\beta }{K^2}, \frac{K^2  \ell c}{8\pi  n_1n_2m \alpha}\Big)  U_3\Big(\frac{n_1 n_2 m \alpha }{\ell ^2 N_3}\Big) \Big(\frac{n_1 n_2 m \alpha }{\ell ^2 N_3}\Big)^{\frac14}  \\
&+O(K^{-1000}) + O\Big( \frac{1}{K^{3-\ep}} \delta_{|\ell| \asymp \frac{\sqrt{n_1n_2m\alpha} }{\sqrt{N_3} }} \Big),
\end{align*}
with the understanding that the main term vanishes if $\ell =0$. The delta function $\delta_{\text{P}}$ equals 1 if the statement P holds and 0 otherwise.
\end{lemma}

\proof

Suppose first that $\ell=0$. Then integrate by parts $j$ times the integral $I$ given in (\ref{od2}). Here we repeatedly integrate $e(\frac{2\sqrt{x n_1n_2N_3 m\alpha}}{c})$ after substituting $y=\sqrt{x}$ and differentiate the rest of the integrand. Using (\ref{wbound}) and (\ref{crange}) we have that 
\begin{align*}
I\ll K^\ep \Big( \frac{ \sqrt{n_1n_2N_3 m \alpha}}{c} \Big)^{-j}   \ll K^\ep \Big(\frac{1}{K^{2-\ep}}\Big)^j.
\end{align*}
Taking $j$ large enough, this is $O(K^{-1000})$.

Now suppose that $\ell \neq 0$ and $|\ell|> K^{100}$. We integrate by parts $j$ times the integral $I$ given in (\ref{od2}). This time we repeatedly integrate $e(\frac{-\ell N_3 x}{c}) $ and differentiate the rest of the integrand. We again get that $I\ll K^{-1000}$.

Henceforth assume $\ell \neq 0$ and $|\ell |\le K^{100}$. Making the substitution
\begin{align*}
y=\frac{\ell ^2 N_3}{n_1 n_2  m \alpha}x,
\end{align*}
the integral is
\begin{multline*}
I= \int_{-\infty}^{\infty} e\Big(\frac{n_1n_2m\alpha }{\ell c}(2\sqrt{y}-y)\Big)  \Psi\Big(\frac{n_1\beta_1}{K},\frac{n_2\beta_2}{K},
\frac{y n_1n_2m\alpha}{\ell^2 K}, \frac{m\beta }{K^2}, \frac{K^2  \ell c}{8\pi  \sqrt{y} n_1n_2 m \alpha}\Big)\\
 \Big(\frac{n_1 n_2 m \alpha}{\ell ^2 N_3} \Big)^{\frac14} U_3\Big(\frac{n_1 n_2 m \alpha }{\ell ^2 N_3}y\Big)  \frac{dy}{y^{\frac{3}{4}}}. 
\end{multline*}
Define $U_0(x)=x^{\frac14} U_3(x)$, so that 
\begin{align}
I\label{idef} = \int_{-\infty}^{\infty} e\Big(\frac{n_1n_2m\alpha }{\ell c}(2\sqrt{y}-y)\Big)  \Psi\Big(\frac{n_1\beta_1}{K},\frac{n_2\beta_2}{K},
\frac{y n_1n_2m\alpha}{\ell^2 K}, \frac{m\beta }{K^2}, \frac{K^2  \ell c}{8\pi  \sqrt{y} n_1n_2 m \alpha}\Big)  U_0\Big(\frac{n_1 n_2 m \alpha }{\ell ^2 N_3}y\Big)  \frac{dy}{y}. 
\end{align}
The stationary point occurs at $y=1$. Let $\Omega(y)$ be a smooth function with bounded derivatives which is equal to 1 on $(1/2, 3/2)$ and 0 on $(-\infty,1/4)\cup(2,\infty)$. We write
\begin{align*}
I=I_1+I_2,
\end{align*}
where $I_1$ is defined as in (\ref{idef}) except that its integrand has an extra factor $1-\Omega(y)$, and $I_2$ is defined as in (\ref{idef}) except that its integrand has an extra factor $\Omega(y)$.

We first show that
\begin{align*}
I_1\ll K^{-1000}.
\end{align*}
For this we will use \cite[lemma 8.1]{bky} with
\begin{align*}
&h(y)= \frac{2\pi n_1n_2m\alpha }{\ell c}(2\sqrt{y}-y), \\
&w(y) = \frac{1-\Omega(y)}{y} \Psi \Big(\frac{n_1\beta_1}{K},\frac{n_2\beta_2}{K},
\frac{y n_1n_2m\alpha}{\ell^2 K}, \frac{m\beta }{K^2}, \frac{K^2  \ell c}{8\pi  \sqrt{y} n_1n_2 m \alpha}\Big)  U_0\Big(\frac{n_1 n_2 m \alpha }{\ell ^2 N_3}y\Big).
\end{align*}
The parameters in this lemma are 
\begin{align*}
R= \frac{n_1 n_2 m \alpha}{|\ell | c} \Big( \Big(\frac{n_1 n_2 m\alpha }{\ell ^2 N_3}\Big)^\half +1\Big), \ \ U=Q= \frac{\ell ^2 N_3}{n_1 n_2 m \alpha},  \ \ Y= \frac{n_1 n_2 m\alpha }{|\ell | c} \Big( \frac{\ell ^2 N_3}{n_1 n_2 m \alpha}\Big)^\half.
\end{align*}
This is because
\begin{align*}
 h'(y)=\frac{n_1n_2 m \alpha}{\ell  c}\Big( \frac{1}{\sqrt{y}} -1 \Big), \ \ h^{(j)}(y) \ll \frac{n_1n_2m\alpha}{|\ell | c} y^\half y^{-j}
\end{align*}
for $j\ge 2$, and we can assume that $|y-1|\gg 1$ by the support of $1-\Omega(y)$ ,and that $y\asymp \frac{\ell ^2 N_3}{n_1 n_2 m \alpha }$ by the support of $U_0$. Further, by (\ref{wbound}), we have
\begin{align*}
w^{(j)}(y)\ll \frac{K^\ep}{y^{j+1}}.
\end{align*}
 We don't need to specify the remaining parameters $\alpha, \beta, X$ given in \cite[lemma 8.1]{bky}, apart from noting that they are bounded by some power of $K$. The result of the lemma is
\begin{align*}
I_1\ll (\beta-\alpha)X( (QR/\sqrt{Y} )^{-A}+(RU)^{-A})
\end{align*}
for any $A\ge 0$. Thus it suffices to show that $QR/\sqrt{Y}>K^\ep$ and $RU>K^\ep$, and then to take $A$ large enough.

{\it Case 1.} Suppose that $ \frac{\ell ^2 N_3}{n_1 n_2 m \alpha}\ge 1$. Then $R\gg \frac{n_1 n_2 m \alpha}{|\ell | c}$ and so
\begin{align*}
&QR/\sqrt{Y} \gg \Big(\frac{n_1 n_2 m \alpha}{|\ell | c}\Big)^\half \Big(  \frac{\ell ^2 N_3}{n_1 n_2 m \alpha} \Big)^\frac34 \gg  \Big(\frac{n_1 n_2 m \alpha}{|\ell | c}\Big)^\half \Big(  \frac{\ell ^2 N_3}{n_1 n_2 m \alpha} \Big)^\frac14  = \Big(\frac{
\sqrt{n_1 n_2 N_3 m \alpha}}{c} \Big)^\half,\\
&RU \gg \frac{n_1 n_2 m \alpha}{|\ell | c} \Big(  \frac{\ell ^2 N_3}{n_1 n_2 m \alpha} \Big) \gg \frac{n_1 n_2 m \alpha}{|\ell | c} \Big(  \frac{\ell ^2 N_3}{n_1 n_2 m \alpha} \Big) ^\frac12 =  \frac{\sqrt{n_1 n_2 N_3 m \alpha}}{c}.
\end{align*}
By (\ref{crange}), we have $\frac{\sqrt{n_1 n_2 N_3 m \alpha}}{c} \gg K^{2-\ep}$. 

{\it Case 2.} Suppose that $ \frac{\ell ^2 N_3}{n_1 n_2 m \alpha}< 1$. Then $R\gg \frac{n_1n_2 m \alpha}{|\ell | c} (\frac{n_1n_2 m \alpha}{\ell ^2 N_3})^\half$ and so
\begin{align*}
&QR/\sqrt{Y} \gg \Big(\frac{n_1 n_2 m \alpha}{|\ell | c}\Big)^\half \Big(  \frac{\ell ^2 N_3}{n_1 n_2 m \alpha} \Big)^\frac14  = \Big(\frac{
\sqrt{n_1 n_2 N_3 m \alpha}}{c} \Big)^\half,\\
&RU \gg  \frac{n_1 n_2 m \alpha}{|\ell | c} \Big(  \frac{\ell ^2 N_3}{n_1 n_2 m \alpha} \Big) ^\frac12 =  \frac{\sqrt{n_1 n_2 N_3 m \alpha}}{c},
\end{align*}
and the conclusion is the same.

\smallskip

Now consider $I_2$. We have
\begin{multline*}
I_2 =  \int_{-\infty}^{\infty} e\Big(\frac{n_1n_2m\alpha }{\ell c}(2\sqrt{y}-y)\Big)  U\Big(\frac{n_1\beta_1}{N_1},\frac{n_2\beta_2 }{N_2},
\frac{y n_1n_2m\alpha}{\ell^2 N_3}, \frac{m\beta }{K^2}, \frac{K^2  \ell c}{8\pi  \sqrt{y} n_1n_2 m \alpha}\Big)  U_0\Big(\frac{n_1 n_2 m \alpha }{\ell ^2 N_3}y\Big)  \Omega(y) \frac{dy}{y}. 
\end{multline*}
We apply \cite[Proposition 8.2]{bky}, with
\begin{align*}
&h(y)= \frac{2\pi n_1n_2m\alpha }{\ell c}(2\sqrt{y}-y), \\
&w(y) = \frac{\Omega(y)}{y} \Psi \Big(\frac{n_1\beta_1}{K},\frac{n_2\beta_2}{K},
\frac{y n_1n_2m\alpha}{\ell^2 K}, \frac{m\beta }{K^2}, \frac{K^2  \ell c}{8\pi  \sqrt{y} n_1n_2 m \alpha}\Big)  U_0\Big(\frac{n_1 n_2 m \alpha }{\ell ^2 N_3}y\Big),\\
&X=K^\ep, V=V_1=Q=1, Y= \frac{n_1 n_2 m\alpha}{|\ell |c}\asymp \frac{ \sqrt{n_1 n_2 N_3 m \alpha}}{c}.
\end{align*}
The approximation to $Y$ is given by $\ell^2 \asymp \frac{n_1 n_2 m \alpha y}{N_3}$ and $y\asymp 1$ by the support of $U_0$ and $\Omega$ respectively. By (\ref{crange}), we have that $Y\gg K^{2-\ep}$. Thus the conditions \cite[line (8.7)]{bky} are satisfied for $\delta=1/5$ say, and we get (we have a factor of $e(-1/8)$ instead of $e(1/8)$ because the second derivative of $h$ is negative)
\begin{align*}
I_2 =  \sqrt{\frac{2\ell c}{ n_1 n_2 m\alpha }}  e\Big(\frac{n_1 n_2 m \alpha}{\ell c}-\frac{\pi}{8}\Big) \Psi \Big(\frac{n_1\beta_1}{K},\frac{n_2\beta_2}{K},
\frac{ n_1n_2m\alpha}{\ell^2 K}, \frac{m\beta }{K^2}, \frac{K^2  \ell c}{8\pi  n_1n_2 m \alpha}\Big)  U_0\Big(\frac{n_1 n_2 m \alpha }{\ell ^2 N_3}\Big) +\text{error},
\end{align*}
where
\begin{align*}
\text{error}=O\Big( \frac{1}{ \sqrt{Y}}  \sum_{1\le n\le 10^6} p_n(1) + K^{-100}\Big),
\end{align*}
where
\begin{align*}
p_n(1)=\frac{1}{n!} \frac{|G^{(2n)}(1)| }{Y^n} , \ \ \ G(t)=w(t)e(H(t)), \ \ \ H(t) = h(t)-h(1)-\tfrac12 h''(1)(t-1)^2.
\end{align*}
Note that $H(1)=H'(1)=H''(1)=0$, and so $G^{(2n)}(1)\ll Y^{^{\lfloor  \frac{2n}{3} \rfloor}} \ll Y^{n-1}$. Thus 
\begin{align*}
\text{error}= O\Big(Y^{-3/2} \delta_{|\ell| \asymp \frac{\sqrt{n_1 n_2 m\alpha} }{\sqrt{N_3} }}\Big) = O\Big( \frac{1}{K^{3-\ep}} \delta_{|\ell| \asymp \frac{\sqrt{n_1 n_2 m\alpha} }{\sqrt{N_3} }}\Big).
\end{align*}
\endproof
Now we are ready to return to (\ref{od2}). We evaluate the $a$-sum there as
\begin{align*}
\sum_{a \bmod c} S(a,n_1n_2m\alpha, c)e\Big(\frac{a\ell}{c}\Big) =
\begin{cases}
c e(\frac{-n_1n_2m\alpha\overline{\ell}}{c}) &\text{ if } (\ell,c)=1\\
0 &\text{ otherwise},
\end{cases}
\end{align*}
and then apply Lemma \ref{stationary} for the integral. The error term of this lemma contributes, using (\ref{condition}) and (\ref{crange}), at most
\begin{align*}
\frac{1}{K^{3-\ep}} \sum_{\substack{n_1\asymp N_1 ,n_2\asymp N_2 \\ m<K^{2+\ep}/\beta  \\ c< \sqrt{N_1N_2N_3m \alpha}/K^{2-\ep} \\ |\ell|\asymp  \sqrt{n_1 n_2 m\alpha / N_3}  }}   \frac{ N_3^\frac14}{ c^\half (n_1 n_2 m)^\frac34 }  \ll K^\ep.
\end{align*}
Thus we only need to consider the contribution of the main term. It suffices to prove (we only treat the terms with $\ell >0$)
\begin{multline}
\label{reqod3} OD_3 := \sum_{\substack{n_1,n_2,m,c,\ell  \ge 1} } e\Big(\frac{-n_1 n_2 m \alpha \overline{\ell }}{c}\Big) e\Big(\frac{n_1 n_2 m \alpha}{\ell c} \Big)   \frac{\tau(m )}{ n_1 n_2  m } \\
 \Psi \Big(\frac{n_1\beta_1}{K},\frac{n_2\beta_2}{K},
\frac{ n_1n_2m\alpha}{\ell^2 K}, \frac{m\beta }{K^2}, \frac{K^2  \ell c}{8\pi  n_1n_2 m \alpha}\Big)  U_1\Big(\frac{n_1  }{N_1}\Big) U_2\Big(\frac{n_2 }{N_2}\Big) U_3\Big(\frac{n_1 n_2 m \alpha }{\ell ^2 N_3}\Big) \ll \alpha K^{2\theta+\ep},
\end{multline}
where it is understood that the sum is restricted to $(\ell ,c)=1$. By the reciprocity relation for exponentials, we have
\begin{multline*}
 OD_3 = \sum_{\substack{n_1,n_2,m,c,\ell  \ge 1} } e\Big( \frac{n_1 n_2 m \alpha \overline{c }}{\ell} \Big)  \frac{\tau(m )}{ n_1 n_2  m } \\
 \Psi \Big(\frac{n_1\beta_1}{K},\frac{n_2\beta_2}{K},
\frac{ n_1n_2m\alpha}{\ell^2 K}, \frac{m\beta }{K^2}, \frac{K^2  \ell c}{8\pi  n_1n_2 m \alpha}\Big)  U_1\Big(\frac{n_1  }{N_1}\Big) U_2\Big(\frac{n_2 }{N_2}\Big) U_3\Big(\frac{n_1 n_2 m \alpha }{\ell ^2 N_3}\Big) .
\end{multline*}

\section{Voronoi summation and fake main terms}

The next goal is to perform Voronoi summation on $m$ but we cannot do so immediately because in the exponential $e(\frac{n_1 n_2 m \alpha \overline{c }}{\ell}) $, the integers $n_1n_2\alpha$ and $\ell $ may not be coprime. We first prepare by eliminating any common factors. Re-ordering the sum $OD_3$ by $b_1= (n_1,\ell )$, and replacing $n_1$ by $b_1 n_1$ and $\ell $ by $b_1 \ell $, we have
\begin{multline*}
OD_3 = \sum_{\substack{n_1,n_2,m,c,\ell  \ge 1\\ b_1\ge1 \\ (n_1,\ell)=1}} e\Big( \frac{n_1 n_2 m \alpha \overline{c }}{\ell} \Big)   \frac{\tau(m )}{ b_1 n_1 n_2  m } \\
 \Psi \Big(\frac{b_1n_1\beta_1}{K},\frac{n_2\beta_2}{K}, \frac{ n_1n_2m\alpha}{b_1 \ell^2 K}, \frac{m\beta }{K^2}, \frac{K^2  \ell c}{8\pi  n_1n_2 m \alpha}\Big)  U_1\Big(\frac{b_1n_1  }{N_1}\Big) U_2\Big(\frac{b_2 n_2 }{N_2}\Big) U_3\Big(\frac{n_1 n_2 m \alpha }{b_1 \ell ^2 N_3}\Big).
\end{multline*}
Next we re-order the sum by $b_2= (n_2,\ell )$, and replace $n_2$ by $b_2 n_2$ and $\ell $ by $b_2 \ell $, then re-order the result by $b_3= (\alpha,\ell )$, and replace $\ell$ by $b_3 \ell$ and $\alpha$ by $b_3 \alpha$. In this way, the conditions (\ref{condition}) become
\begin{align}
\label{condition2} N_3\ge N_2, \ \ N_1N_2<\frac{K^{2+\ep}}{b_3 \alpha}, \ \ \ \beta\ge b_3\alpha,
\end{align}
and
\begin{multline}
\label{gcd} OD_3 = \sum_{\substack{n_1,n_2,m,c,\ell  \ge 1\\ b_1,b_2\ge1,b_3|\alpha \\ (n_1n_2\alpha,\ell)=1}}  e\Big( \frac{n_1 n_2 m \alpha \overline{c }}{\ell} \Big)   \frac{\tau(m )}{ b_1 b_2  n_1 n_2  m } \\
 \Psi \Big(\frac{b_1n_1\beta_1}{K},\frac{b_2 n_2\beta_2}{K}, \frac{ n_1n_2m\alpha}{b_1b_2b_3  \ell^2 K}, \frac{m\beta }{K^2}, \frac{K^2  \ell c}{8\pi  n_1n_2 m \alpha}\Big)  U_1\Big(\frac{b_1 n_1  }{N_1}\Big) U_2\Big(\frac{n_2 b_2 }{N_2}\Big) U_3\Big(\frac{n_1 n_2 m \alpha }{b_1 b_2 b_3\ell ^2 N_3}\Big),
\end{multline}
for which the required bound (\ref{reqod3}) becomes
\begin{align*}
OD_3\ll b_3 \alpha K^{2\theta+\ep}.
\end{align*}

\smallskip 

Working in dyadic intervals of $m$ by taking a partition of unity, it suffices to show that
\begin{multline*}
OD_3 = \sum_{j} \sum_{\substack{n_1,n_2,m,c,\ell  \ge 1\\ b_1,b_2\ge1,b_3|\alpha \\ (n_1n_2\alpha,\ell)=1}}  e\Big( \frac{n_1 n_2 m \alpha \overline{c }}{\ell} \Big)   \frac{\tau(m )}{ b_1 b_2  n_1 n_2  m }  U_1\Big(\frac{b_1 n_1  }{N_1}\Big) U_2\Big(\frac{n_2 b_2 }{N_2}\Big) U_3\Big(\frac{n_1 n_2 m \alpha }{b_1 b_2 b_3\ell ^2 N_3}\Big) U_{4,j}\Big(\frac{m  }{M_j}\Big)\\
 \Psi \Big(\frac{b_1n_1\beta_1}{K},\frac{b_2 n_2\beta_2}{K}, \frac{ n_1n_2m\alpha}{b_1b_2b_3  \ell^2 K}, \frac{m\beta }{K^2}, \frac{K^2  \ell c}{8\pi  n_1n_2 m \alpha}\Big) \ll b_3 \alpha K^{\theta+\ep}
\end{multline*}
for some smooth functions $U_{4,j}$ compactly supported on $(\half,\frac52)$, for $j\ll \log K$, and $M_j\asymp 2^j$. We now apply the Voronoi summation formula (Lemma \ref{vorsum}) to the $m$ sum, getting
\begin{align}
\label{od3def} OD_3:= FM + OD_4,
\end{align}
where the ``fake main term'' is
\begin{multline*}
FM :=  \int_{-\infty}^{\infty} \Big( \log \frac{x}{\ell ^2} +2\gamma \Big) \sum_{j} U_{4,j}\Big(\frac{x  }{M_j}\Big)  \sum_{\substack{n_1,n_2,c,\ell  \ge 1\\ b_1,b_2\ge1,b_3|\alpha \\ (n_1n_2\alpha,\ell)=1}}    \frac{1}{\ell  b_1 b_2  n_1 n_2   } \\
  \Psi \Big(\frac{b_1n_1\beta_1}{K},\frac{b_2 n_2\beta_2}{K}, \frac{ n_1n_2x\alpha}{b_1b_2b_3  \ell^2 K}, \frac{x\beta }{K^2}, \frac{K^2  \ell c}{8\pi  n_1n_2 x \alpha}\Big) U_1\Big(\frac{b_1 n_1  }{N_1}\Big) U_2\Big(\frac{n_2 b_2 }{N_2}\Big) U_3\Big(\frac{n_1 n_2 x \alpha }{b_1 b_2 b_3\ell ^2 N_3}\Big) \frac{dx}{x},
\end{multline*}
and $OD_4$ is given in the next section. In the sum $FM$, we may re-patch the partition of unity and reverse the steps which led to (\ref{gcd}), getting that
\begin{multline*}
FM =   \int_{-\infty}^{\infty} \Big( \log \frac{x}{\ell ^2} +2\gamma \Big)  \sum_{\substack{n_1,n_2,c,\ell  \ge 1} }  \frac{1}{ \ell n_1 n_2   } \\
 \Psi \Big(\frac{n_1\beta_1}{K},\frac{n_2\beta_2}{K},
\frac{ n_1n_2x\alpha}{\ell^2 K}, \frac{x\beta }{K^2}, \frac{K^2  \ell c}{8\pi  n_1n_2  x \alpha}\Big)  U_1\Big(\frac{n_1  }{N_1}\Big) U_2\Big(\frac{n_2 }{N_2}\Big) U_3\Big(\frac{n_1 n_2 x \alpha }{\ell ^2 N_3}\Big) \frac{dx}{x}.
\end{multline*}
The trivial bound for $M$ is $O(K^{\half+\ep})$, from the length of the $c$-sum given by (\ref{crange}). It seems like we cannot do better because there are no exponentials or other harmonics present which may produce further cancellation (hence the name ``fake main term''). However we can exploit our judicious choice of weight function in the approximate functional equation, as follows.

Making the substitution $y=\frac{x n_1 n_2 \alpha}{\ell ^2 N_3}$, we have
\begin{multline*}
FM =   \int_{-\infty}^{\infty}  \sum_{\substack{n_1,n_2,c,\ell  \ge 1} }  \Big( \log \frac{yN_3}{n_1n_2\alpha} +2\gamma \Big)   \frac{1}{ \ell n_1 n_2   } \\
 \Psi \Big(\frac{n_1\beta_1}{K},\frac{n_2\beta_2}{K},
\frac{ y N_3}{K}, \frac{y \ell^2 N_3 \beta }{K^2 n_1n_2\alpha}, \frac{K^2  c}{8\pi  \ell N_3 y}\Big)  U_1\Big(\frac{n_1  }{N_1}\Big) U_2\Big(\frac{n_2 }{N_2}\Big) U_3(y) \frac{dy}{y}.
\end{multline*}
It suffices to show that
\begin{align*}
 FM':= \sum_{\substack{ c,\ell  \ge 1} }  \frac{1}{  \ell  } 
 \Psi \Big(\frac{n_1\beta_1}{K},\frac{n_2\beta_2}{K},
\frac{ y N_3}{K}, \frac{y \ell^2 N_3 \beta }{K^2 n_1n_2\alpha}, \frac{K^2  c}{8\pi  \ell N_3 y}\Big)  \ll K^\ep
\end{align*}
for any $n_1\asymp N_1, n_2\asymp N_2, y\asymp 1$. Using (\ref{wdef}) and Mellin inversion, we have
\begin{multline*}
FM'= \frac{1}{(2\pi i)^5} \int_{(1+\ep)} \int_{(\half+\ep)} \int_{(\ep)} \int_{(\ep)}  \int_{(\ep)}   \frac{e^{s_1^2+s_2^2+s_3^2} \zeta(1+2s_4) \mathcal{G}(s_4) }{ 2^{s_1+s_2+s_3} 4^{s_4}  }  \zeta(1+2s_4-w) \zeta(w)  \\ 
\Big(\frac{K}{n_1\beta_1}\Big)^{s_1}\Big(\frac{K}{n_2\beta_2}\Big)^{s_2}  \Big(\frac{K}{yN_3}\Big)^{s_3}  \Big( \frac{K^2 n_1n_2\alpha}{y  N_3 \beta }\Big)^{s_4}  \Big(\frac{8\pi   N_3 y}{K^2  }\Big)^w   
\tilde{\hbar}_{s_1+2s_2+2s_3}(w) \frac{ds_1}{s_1}\frac{ds_2}{s_2}\frac{ds_3}{s_3}\frac{ds_4}{s_4} dw.
\end{multline*} 
Here $\zeta(1+2s_4-w)$ comes from the $\ell $-sum and $\zeta(w)$ comes from the $c$-sum. We must initially keep the lines of integration at $\Re(w)=1+\ep$ and $\Re(s_4)=\half+\ep$ in order to stay in the region of absolute convergence. The goal is to move all the lines of integration to $(\ep)$, and this would prove the claim.
We first move the $w$-integral to $\Re(w)=\ep$. This crosses a simple pole at $w=1$, with residue
\begin{multline*}
FM'':= \frac{1}{(2\pi i)^4 } \int_{(\half+\ep)} \int_{(\ep)} \int_{(\ep)}  \int_{(\ep)}   \frac{e^{s_1^2+s_2^2+s_3^2} \zeta(1+2s_4) \mathcal{G}(s_4) }{ 2^{s_1+s_2+s_3} 4^{s_4}  }  \zeta(2s_4)  \\ 
\Big(\frac{K}{n_1\beta_1}\Big)^{s_1}\Big(\frac{K}{n_2\beta_2}\Big)^{s_2}  \Big(\frac{K}{yN_3}\Big)^{s_3}  \Big( \frac{K^2 n_1n_2\alpha}{y  N_3 \beta }\Big)^{s_4}  \Big(\frac{8\pi   N_3 y}{K^2  }\Big)   
\tilde{\hbar}_{s_1+2s_2+2s_3}(1) \frac{ds_1}{s_1}\frac{ds_2}{s_2}\frac{ds_3}{s_3}\frac{ds_4}{s_4} .
\end{multline*} 
On the shifted integral at $\Re(w)=\ep$, which is not displayed, we may move the $s_4$ integral to $\Re(s_4)=\ep$ and then estimate (this does not cross any pole of $\zeta(1+2s_4-w)$ so this straight forward). Thus the shifted integral is $O(K^{\ep})$ and we are left to estimate $FM''$. In the integral $FM''$, we move the line of integration to $\Re(s_4)=\ep$. This does not cross any poles because the simple pole of $\zeta(2s_4)$ at $s_4=\half$ is cancelled out by the zero at $s_4=\half$ of $\mathcal{G}(s_4)$. See the definition (\ref{gdef}). Thus $FM''$ is $O(K^{-1+\ep})$.

\section{Second application of reciprocity}

We now return to (\ref{od3def}) and give the definition of $OD_4$ corresponding to the sum on the right hand side of (\ref{vorstate}). We have ($r$ is the dual variable)
\begin{multline*}
OD_4 := \sum_{j}\sum_\pm \frac{1}{2\pi i} \int_{(A)} \int_0^\infty\sum_{\substack{n_1,n_2,r,c,\ell  \ge 1\\ b_1,b_2\ge1,b_3|\alpha \\ (c,\ell)=1}}  e\Big( \frac{\pm r c \overline{n_1 n_2  \alpha} }{\ell} \Big)   \frac{\tau(r )}{ \ell b_1 b_2  n_1 n_2   }   \Big(\frac{M_j r x}{\ell ^2}\Big)^{-w} H_1^{\pm}(w) U_1\Big(\frac{b_1 n_1  }{N_1}\Big)\\
 U_2\Big(\frac{n_2 b_2 }{N_2}\Big)  U_3\Big(\frac{n_1 n_2 x M_j \alpha }{b_1 b_2 b_3\ell ^2 N_3}\Big)    U_{4,j}(x) 
 \Psi \Big(\frac{b_1n_1\beta_1}{K},\frac{b_2 n_2\beta_2}{K}, \frac{ n_1n_2xM_j\alpha}{b_1b_2b_3  \ell^2 K}, \frac{xM_j\beta }{K^2}, \frac{K^2  \ell c   }{8\pi    n_1   n_2 x M_j \alpha}\Big) \frac{dx}{x} dw,
\end{multline*}
where it is understood that the sum is restricted to $(n_1n_2\alpha,\ell)=1$ and we need $OD_4\ll b_3 \alpha K^{2\theta+\ep}$.

We first simplify the notation a bit (we did not do this earlier because we needed the exact form of the weight functions in order to deal with fake main terms). First, we observe that since there are $O(K^\ep)$ dyadic intervals, it is enough to consider any one smooth function $U_{4,j}=U_4$ and $M_j=M$. From the fourth component of $\Psi$ and the assumption $\beta\ge b_3 \alpha$ from (\ref{condition}), we can assume
\begin{align*}
M < \frac{K^{2+\ep}}{b_3\alpha}.
\end{align*}
 We can also consider the sum in dyadic intervals $r\asymp R$ by inserting a smooth bump function $U_5(\frac{r}{R})$, where $U_5$ is supported on $(\half,\frac52)$. We can assume that
\begin{align*}
 R<\frac{K^\ep  \ell^2}{ M}
\end{align*}
because the contribution of $r \ge \frac{K^\ep  \ell^2}{ M}$ is $O(K^{-100})$ say. This can be seen by moving the $w$-integral in $OD_4$ far to the right (taking $A$ large). By repeatedly integrating by parts the $x$-integral, we may restrict the $w$-integral to $|\Im w|<K^\ep$ (the real part is already fixed at $A$). Doing so, we may absorb $r^{-w}$ and $(\ell^2)^{-w}$ into $U_5$ and $U_3$ respectively. Similarly we may expand the function $\Psi$ using (\ref{wdef}), truncate the integrals there to $|\Im s_1|,|\Im s_2|,|\Im s_3|,|\Im s_4|<K^\ep$ (with $\Re(s_i)$ fixed of course) and absorb part of this function into the bump functions $U_1$, $U_2$, $U_3$, $U_4$. Thus it suffices to prove (we do not seek cancellation in the sum over $b_1,b_2,b_3$)
\begin{align*}
 \sum_{\substack{n_1,n_2,r,c,\ell  \ge 1\\  \\ (c,\ell)=1}}  e\Big( \frac{\pm r c \overline{n_1 n_2  \alpha} }{\ell} \Big)   \frac{\tau(r )}{ \ell  n_1 n_2   } U_1\Big(\frac{b_1 n_1  }{N_1}\Big)  U_2\Big(\frac{n_2 b_2 }{N_2}\Big)
  U_3\Big(\frac{n_1 n_2  M \alpha }{b_1 b_2 b_3 \ell ^2 N_3}\Big)  U_5\Big( \frac{r}{R}\Big) \hbar_s \Big(\frac{K^2  \ell c}{8\pi x n_1n_2  M \alpha}\Big) \\
  \ll b_3 \alpha K^{\theta+\ep}
\end{align*}
for any $b_1,b_2,b_3\ge 1$, $x\asymp 1$, $|s|<K^\ep$ and any compactly supported functions $U_j$ with $j$-th derivative bounded by $(K^\ep)^j$. We simplify the notation a bit more. We suppress the factor $8\pi x$ in $\hbar_s$, rename $b_1b_2$ to $a$, $b_3$ to $b$, $M\alpha$ to $M$, $N_1/b_1$ to $N_1$ and $N_2/b_2$ to $N_2$. Thus it suffices to prove
\begin{align}
 \label{b4recip} \sum_{\substack{n_1,n_2,r,c,\ell  \ge 1\\  \\ (c,\ell)=1}}  e\Big( \frac{\pm r c \overline{n_1 n_2  \alpha} }{\ell} \Big)   \frac{\tau(r )}{ \ell  n_1 n_2   } U_1\Big(\frac{ n_1  }{N_1}\Big)  U_2\Big(\frac{ n_2 }{N_2}\Big)
  U_3\Big(\frac{n_1 n_2  M  }{ab \ell ^2 N_3}\Big)  U_5\Big( \frac{r}{R}\Big) \hbar_s \Big(\frac{K^2  \ell c}{ n_1n_2  M}\Big)   \ll b \alpha K^{\theta+\ep},
\end{align}
for any integers $a,b,\alpha$ and 
\begin{align}
\label{newcond} N_1,N_2,N_3 < K^{1+\ep}, \ \ N_1N_2<\frac{K^{2+\ep}}{\alpha a b}, \ \ M<\frac{K^{2+\ep}}{b},  \ \ 
R<\frac{K^\ep  \ell^2 \alpha}{ M } \asymp \frac{K^\ep  N_1 N_2 \alpha}{ ab N_3}, 
 \ \ N_3\ge N_2.
\end{align}
The approximation $\frac{K^\ep  \ell^2 \alpha}{ M } \asymp \frac{K^\ep  N_1 N_2 \alpha}{ ab N_3} $ follows from the support of $U_3$. This updates (\ref{condition2}).

\bigskip 

Now using reciprocity for exponentials, we have
\begin{align*}
e\Big( \frac{ \pm r c \overline{n_1 n_2  \alpha} }{\ell} \Big)  =e\Big( \frac{ \mp r c \overline{\ell} }{n_1 n_2  \alpha} \Big) e\Big( \frac{ \pm r c  }{\ell n_1 n_2  \alpha}  \Big)  = e\Big( \frac{ \mp r c \overline{\ell} }{n_1 n_2  \alpha} \Big)  \Big( 1+ O\Big( \frac{ \pm r c  }{\ell n_1 n_2  \alpha}\Big) \Big).
\end{align*}
The contribution to (\ref{b4recip}) of this error term is less than
\begin{align*}
\sum_{\substack{n_1\asymp N_1\\ n_2\asymp N_2\\ \ell \asymp \frac{ \sqrt{ N_1 N_2 M}}{\sqrt{a b N_3}} }}  \ \sum_{ r < \frac{K^\ep N_1 N_2 \alpha}{N_3 a b}} \ \sum_{c< \frac{ N_1N_2 M  }{\ell K^{2-\ep}}} \frac{K^\ep}{\ell n_1n_2} \cdot \frac{ r c  }{\ell n_1 n_2  \alpha}
\ll \frac{N_1^\frac32 N_2^\frac32 M^\half \alpha}{K^{4-\ep} N_3^\half}\ll \frac{1}{K^{\half-\ep}},
\end{align*}
by (\ref{newcond}).
So in (\ref{b4recip}) we can replace the exponential with $e( \frac{ \mp r c \overline{\ell} }{n_1 n_2  \alpha})$ and detect the condition $(\ell ,c)=1$ using the M\"{o}bius function:
\begin{align}
\sum_{l |(\ell,c)} \mu(l) = 
\begin{cases}
1 &\text{ if } (\ell,c)=1\\
0 &\text{otherwise}.
\end{cases}
\end{align}
Thus replacing $\ell$ by $\ell l$ and $c$ by $cl$, 
 it suffices to prove
\begin{align*}
 \sum_{\substack{n_1,n_2,r,c,\ell  \ge 1\\  l\ge 1 }  }e\Big( \frac{ \mp r c \overline{\ell} }{n_1 n_2  \alpha} \Big)    \frac{\mu(l) \tau(r )}{ l \ell  n_1 n_2   } U_1\Big(\frac{ n_1  }{N_1}\Big)  U_2\Big(\frac{ n_2 }{N_2}\Big)
  U_3\Big(\frac{n_1 n_2  M  }{ab l^2 \ell ^2 N_3}\Big)  U_5\Big( \frac{r}{R}\Big) \hbar_s \Big(\frac{K^2 l^2 \ell c}{ n_1n_2  M}\Big) 
  \nonumber  \ll b \alpha K^{\theta+\ep}.
\end{align*}
We do not seek cancellation over the $l$-sum, so it suffices to prove
\begin{align*}
OD_5:= \sum_{\substack{n_1,n_2,r,c,\ell  \ge 1 }  } e\Big( \frac{ \mp r c \overline{\ell} }{n_1 n_2  \alpha} \Big)   \frac{ \tau(r )}{  \ell  n_1 n_2   } U_1\Big(\frac{ n_1  }{N_1}\Big)  U_2\Big(\frac{ n_2 }{N_2}\Big)
  U_3\Big(\frac{n_1 n_2  M  }{ab l^2 \ell ^2 N_3}\Big)  U_5\Big( \frac{r}{R}\Big) \hbar_s \Big(\frac{K^2 l^2 \ell c}{ n_1n_2  M}\Big) 
  \nonumber  \ll b \alpha K^{\theta+\ep}
\end{align*}
for any integer $l\ge 1$ and assuming (\ref{newcond}). Also keep in mind that it is understood that the sum is restricted to $(\ell,n_1n_2\alpha)=1$.

\section{Second Poisson summation}

Now we split the $\ell $-sum in $OD_5$ into (primitive) residue classes mod $n_1 n_2 \alpha$ and apply Poisson summation (Lemma \ref{poiss}). Note that $\ell$ is supported in compact interval of size ${\frac{\sqrt{n_1n_2M}}{l\sqrt{abN_3}}}$. The result is that (the dual variable is $d$)
\begin{align}
\label{od5pois} OD_5= &\sum_{\substack{n_1,n_2,r,c \ge 1 }} \frac{S( \mp r c,  d, n_1 n_2  \alpha)}{n_1n_2\alpha}    \frac{\tau(r ) }{   n_1 n_2   } U_1\Big(\frac{ n_1  }{N_1}\Big)  U_2\Big(\frac{ n_2 }{N_2}\Big)
    U_5\Big( \frac{r}{R}\Big)  \\
\nonumber &  \times \Bigg(  \int_{-\infty}^{\infty} U_3\Big(\frac{1}{y}\Big) \hbar_s \Big(\frac{y K^2   c l }{ \sqrt{ab n_1n_2 N_3 M}}\Big)  \frac{dy}{y} \\
\nonumber  & + \sum_{d\neq 0}\frac{1}{2\pi i}\int_{(A)}   \int_0^\infty   H_2(z)  
  \Big(  \frac{-2\pi yd\sqrt{M}}{l \alpha \sqrt{abn_1n_2N_3}} \Big)^{-z}
U_3\Big(\frac{1}{y}\Big) \hbar_s \Big(\frac{y K^2   c l}{ \sqrt{ab n_1n_2 N_3 M}}\Big)  \frac{dy}{y} dz \Bigg).
\end{align}
We first consider the contribution of the second line of (\ref{od5pois}). This is the zero frequency contribution, and it is bounded by
\begin{align*}
\sum_{\substack{n_1\asymp N_1 \\ n_2\asymp N_2}} \sum_{ r < \frac{K^\ep N_1 N_2 \alpha}{N_3 a b}}   \sum_{c< \frac{\sqrt{ab n_1n_2  N_3 M  }}{l K^{2-\ep} } } \frac{K^\ep |S( \mp r c,  0, n_1 n_2  \alpha)|}{n_1^2 n_2^2 \alpha}    \ll 
\frac{\sqrt{N_1N_2M}}{K^{2-\ep}\sqrt{N_3}}\ll \frac{1}{K^{\half-\ep}},
\end{align*}
on using $N_3\ge  N_2$ and that the Ramanujan sum is $O(K^\ep)$ on average. 

Now we consider the contribution of the third line of (\ref{od5pois}), arising from the sum over $d\neq 0$. We consider this sum in dyadic intervals $d\asymp D$ (for simplicity, we restrict to only positive values of $d$) and $c\asymp C$ by inserting smooth bump functions $U_6(\frac{d}{D})$ and $U_7(\frac{c}{C})$ say. We can assume that
\begin{align}
\label{drange} D< \frac{ K^\ep l \alpha \sqrt{ab N_1 N_2 N_3} } { \sqrt{ M } }
\end{align}
because the contribution of $d>\frac{l \alpha \sqrt{ab N_1 N_2 N_3} } { \sqrt{ M } }$ is $O(K^{-100})$ by moving $z$-integral in (\ref{od5pois}) far to the right. Restricting to $|\Im z|<K^\ep$ and $\Re(z)$ fixed, which we may do up to negligible error by repeatedly by parts with respect to $y$, we may absorb $d^{-z},n_1^{z},n_2^{z}$ into the existing weight functions. We can also assume that
\begin{align}
\label{crange2} C< \frac{\sqrt{ab N_1N_2  N_3 M  }}{l K^{2-\ep} }
\end{align}
and absorb the function $\hbar_s$ into $U_7$, by using Mellin inversion and separating variables as above. Thus it suffices to prove
\begin{align*}
\sum_{\substack{n_1,n_2,r,c,d \ge 1 }} \frac{S( \mp r c,  d, n_1 n_2  \alpha)}{n_1n_2\alpha}    \frac{\tau(r ) }{   n_1 n_2   } U_1\Big(\frac{ n_1  }{N_1}\Big)  U_2\Big(\frac{ n_2 }{N_2}\Big)    U_5\Big( \frac{r}{R}\Big) U_6\Big( \frac{d}{D}\Big) U_7\Big( \frac{c}{C}\Big) \ll b\alpha K^{2\theta+\ep}.
\end{align*}
Finally, we need this to be in a form to which we can apply Kuznetsov's formula. To this end we define
\begin{align*}
X:=\frac{N_1N_2\alpha}{\sqrt{RDC}},
\end{align*}
and replace $U_2(\frac{n_2}{N_2})$ with a different bump function
\begin{align*}
Y_1\Big(\frac{  4\pi  X \sqrt{rcd}}{ n_1 n_2 \alpha  }\Big) 
\end{align*}
with properties given below. We can also replace $\frac{\tau(r ) }{   n_1 n_2   }$ with $\frac{\tau(r ) }{   N_1 N_2   }$. Thus it suffices to prove (we do not seek cancellation in the $n_1$ sum)
\begin{align*}
OD_6:= \frac{1}{ N_1 N_2 } \sum_{n_1\asymp N_1}  \Big| \sum_{\substack{n_2 ,r,c,d \ge 1 }}  \frac{S( \pm r c,  d, n_1 n_2 \alpha) }{n_1 n_2 \alpha} 
Y_1\Big(\frac{  4\pi  X \sqrt{rcd}}{ n_1 n_2 \alpha  }\Big)  Y_2\Big(\frac{d }{ D}\Big)  Y_3\Big(\frac{r}{R}\Big) Y_4\Big(\frac{ c }{C}\Big) \Big|  \ll b\alpha K^{2\theta+\ep},
\end{align*}
where $Y_i$ are smooth functions compactly supported on $(\half,\frac{5}{2})$ with $\| Y_j ^{(j)} \|_\infty \ll (K^\ep)^j$ and we assume (\ref{newcond}), (\ref{drange}) and (\ref{crange2}).

\section{Kuznetsov's formula}

The goal now is to prove the required bound for $OD_6$ using Kuznetsov's formula and the spectral large sieve. We consider only the case of positive sign; the negative sign case is similar. By \cite[Theorem 16.5]{iwakow}, we have that 
\begin{align}
\label{kuzback} \sum_{\substack{n_2 \ge 1 }} \frac{ S(  r c,  d, n_1 n_2 \alpha ) }{n_1 n_2 \alpha} 
Y_1\Big(\frac{  4\pi  X \sqrt{rcd}}{ n_1 n_2 \alpha  }\Big) = \text{Maass}  + \text{Eis}+\text{Hol},
\end{align}
where these are the contributions of the Maass cusp forms, Eisenstein series, and holomorphic forms as given in the referenced theorem. For the Maass forms we have
\begin{align*}
\text{Maass} =\sum_{j\ge 1}  \mathcal{M}_{Y_1}(t_j) \frac{\rho_j(rc)\overline{\rho}_j(d)}{\cosh (\pi t_j)},
\end{align*}
where the sum is over an orthonormal basis of Maass cusp forms of level $n_1\alpha$ with Fourier coefficients $\rho_j(n)$ and Laplacian eigenvalue $\frac{1}{4}+t_j^2$, and 
\begin{align*}
 \mathcal{M}_{Y_1}(t) = \frac{\pi i}{2 \sinh( \pi t)} \int_0^\infty (J_{2it}(x)-J_{-2it}(x))  Y_1(xX) \frac{dx}{x}
\end{align*}
By \cite[Lemma 3.6]{butkha}, for example, we have that $ \mathcal{M}_{Y_1}(t) \ll K^{-100}$ if $|t|\ge K^\ep$, so we can restrict the sum Maass to $|t_j|<K^\ep$, in which range $\mathcal{M}_{Y_1}(t_j) \ll  X^{2\theta+\ep}$, by the same Lemma.  
We have $X\ll K^{1-\ep}$ by (\ref{rdc}), so it suffices to prove that
\begin{align*}
\sum_{n_1\asymp N_1} \frac{1}{N_1N_2}  \sum_{|t_j|<K^\ep} \Big| \sum_{\substack{r,d,c\ge 1}} \frac{\rho_j(rc)\overline{\rho}_j(d)}{\cosh(\pi t_j)}  \tau(r)  Y_2\Big(\frac{d }{ D}\Big)  Y_3\Big(\frac{r}{R}\Big) Y_4\Big(\frac{ c }{C}\Big)  \Big| \ll b \alpha K^\ep.
\end{align*}
Now we would like to decompose $\rho_j(rc)$, so that Cauchy-Schwarz and the spectral large sieve may be applied. To do this we need to work with newforms, whose Fourier coefficients are multiplicative. We consult \cite[section 3]{blokha} to see how to choose a basis consisting of lifts of newforms. By \cite[equation (3.10)]{blokha}, and the $\cosh(\pi t_j)^\half$ normalization from the first display of \cite[section 3.2]{blokha}, it suffices to prove that
\begin{align*}
\sum_{n_1\asymp N_1} \frac{1}{N_1N_2}  \sum_{|t_j|<K^\ep}  \frac{K^\ep (uv)^\half}{N_1\alpha} \Big| \sum_{\substack{r,d,c\ge 1\\  u|rc \\ v|d }} \lambda_j\Big(\frac{rc}{u}\Big)\lambda_j\Big(\frac{d}{v}\Big)  \tau(r)  Y_2\Big(\frac{d }{ D}\Big)  Y_3\Big(\frac{r}{R}\Big) Y_4\Big(\frac{ c }{C}\Big)  \Big|\ll b\alpha K^\ep
\end{align*}
for any integers $u,v\ge 1$ and $N_0|n_1\alpha$, where $\lambda_j(n)$ are the Hecke eigenvalues corresponding to newforms of level $N_0$. We now replace $d$ by $dv$ and, proceeding exactly like in steps (\ref{proc1}) to (\ref{proc2}), we can write $u=u_1u_2u_3$ and replace $r$ by $ru_1u_2$ and $c$ by $cu_2u_3$ to see that it suffices to prove
\begin{multline*}
\sum_{n_1\asymp N_1} \frac{1}{N_1N_2}  \sum_{|t_j|<K^\ep}  \frac{K^\ep (u_1u_2u_3v)^\half}{N_1\alpha}\\ \Big| \sum_{\substack{r,d,c\ge 1}} \lambda_j(rcu_2)\lambda_j(d) \mu(u_2) \tau(r u_1 u_2)  Y_2\Big(\frac{d v }{ D}\Big)  Y_3\Big(\frac{r u_1 u_2 }{R}\Big) Y_4\Big(\frac{ cu_2 u_3 }{C}\Big)  \Big|\ll b\alpha K^\ep.
\end{multline*}
To simplify notation, we may replace $r$ by $ru_2$. Thus it suffices to prove
\begin{align*}
\sum_{n_1\asymp N_1} \frac{1}{N_1N_2}  \sum_{|t_j|<K^\ep}  \frac{K^\ep (u_1u_2u_3v)^\half}{N_1\alpha}  \Big| \sum_{\substack{r\asymp R/u_1\\ c\asymp C/u_2u_3\\ d\asymp D/v}} \lambda_j(rc)\lambda_j(d) \gamma_r \gamma_d  \gamma_c   \Big| \ll b\alpha  K^\ep
\end{align*}
for any $\gamma_r,\gamma_c, \gamma_d \ll K^\ep$. By Hecke multiplicativity, we have 
\begin{align*}
\lambda_j(rc)\lambda_j(d) = \sum_{\substack{s|(r,c) \\ (s,N_0)=1}} \mu(s) \lambda_j\Big(\frac{r}{s} \Big)\lambda_j\Big(\frac{c}{s} \Big) \lambda_j(d) = \sum_{\substack{s|(r,c) \\  w|(c/s,d) \\(sw,N_0)=1}} \mu(s) \lambda_j\Big(\frac{r}{s} \Big)\lambda_j\Big(\frac{cd}{sw^2} \Big),
\end{align*}
and so after replacing $r$ by $rs$, $c$ by $csw$, and $d$ by $dw$, it suffices to prove
\begin{align*}
OD_7:= \sum_{\substack{n_1\asymp N_1\\ s,w \le K^{10} }} \frac{1}{N_1N_2}  \sum_{|t_j|<K^\ep}   \frac{K^\ep (u_1u_2u_3v)^\half}{N_1\alpha} \Big| \sum_{\substack{r\asymp R/u_1 s \\ cd \asymp CD/u_2u_3 v sw^2}} \lambda_j(r) \lambda_j(cd) \gamma_{r} \gamma_{cd}     \Big| \ll b\alpha K^\ep,
\end{align*}
for any $\gamma_r, \gamma_{cd}\ll K^\ep$.
By the Cauchy-Schwarz inequality and the spectral large sieve \cite[Theorem 2]{desiwa}, we have that $OD_7$ is bounded by
\begin{align*}
& \sum_{\substack{n_1\asymp N_1\\ s,w \le K^{10} }} \frac{1}{N_1N_2}  \frac{K^\ep (u_1u_2u_3v)^\half}{N_1\alpha}  \Big(  \sum_{|t_j|<K^\ep} \Big| \sum_{\substack{r\asymp R/u_1s }} \lambda_j(r)  \gamma_{r} \Big|^2\Big)^\half   \Big( \sum_{|t_j|<K^\ep} \Big| \sum_{\substack{ cd \asymp CD/u_2u_3vsw^2 }} \lambda_j(cd)  \gamma_{cd}   \Big|^2\Big)^\half \\
&\ll \sum_{\substack{n_1\asymp N_1\\ s,w \le K^{10} }} \frac{1}{N_1N_2}  \frac{K^\ep (u_1u_2u_3v)^\half}{N_1\alpha}   \Big( \Big(N_1\alpha + \frac{R}{u_1s}\Big)\frac{R}{u_1s} \Big)^\half \Big( \Big(N_1\alpha + \frac{CD}{u_2 u_3 v s w^2}\Big)\frac{CD}{u_2 u_3 v s w^2} \Big)^\half.
\end{align*}
Thus it suffices to prove
\begin{align*}
 \frac{K^\ep}{ N_1 N_2 \alpha }   \big( (N_1\alpha + R)R \big)^\half \big((N_1\alpha + CD )CD \big)^\half \ll b\alpha K^\ep.
\end{align*}
By (\ref{newcond}),(\ref{drange}), and (\ref{crange2}), we have
\begin{align}
\label{rdc} (RCD)^\half \ll \frac{\alpha N_1 N_2}{K^{1-\ep}},
\end{align}
so it suffices to prove
\begin{align}
\label{last} \frac{1}{K} (N_1\alpha + R) ^\half  (N_1\alpha + CD )  ^\half \ll b\alpha K^\ep.
\end{align}
We have
\begin{align*}
&\frac{N_1\alpha}{K} \ll K^\ep \alpha,\\ 
& \frac{(N_1\alpha CD )^\half}{K}\ll \frac{N_1\alpha ( a b N_2 N_3 )^\half}{K^2} \ll  \frac{(N_1 N_3  \alpha)^\half  ( a b \alpha N_1 N_2 )^\half}{K^2} \ll K^\ep \alpha^\half,\\
& \frac{(N_1\alpha R )^\half}{K}\ll\frac{N_1\alpha  N_2  ^\half}{K N_3^\half } \ll K^\ep \alpha,
\end{align*}
where in the last inequality we use crucially that $N_3\ge N_2$. This establishes (\ref{last}). 

It remains to consider Eis and Hol in (\ref{kuzback}). These are similarly treated using the large sieve, once we use the multiplicative Fourier coefficients provided explicitly in \cite[section 3]{blokha}.

\bibliographystyle{amsplain}

\bibliography{fifth-moment}

\providecommand{\bysame}{\leavevmode\hbox to3em{\hrulefill}\thinspace}
\providecommand{\MR}{\relax\ifhmode\unskip\space\fi MR }
\providecommand{\MRhref}[2]{%
  \href{http://www.ams.org/mathscinet-getitem?mr=#1}{#2}
}
\providecommand{\href}[2]{#2}
\begin{thebibliography}{10}

\bibitem{bhm}
V.~Blomer, G.~Harcos, and P.~Michel, \emph{Bounds for modular {$L$}-functions
  in the level aspect}, Ann. Sci. \'Ecole Norm. Sup. (4) \textbf{40} (2007),
  no.~5, 697--740.

\bibitem{blokha}
V.~Blomer and R.~Khan, \emph{Twisted moments of $l$-functions and spectral
  reciprocity}, preprint, arXiv:1706.01245.

\bibitem{bky}
V.~Blomer, R.~Khan, and M.~Young, \emph{Distribution of mass of holomorphic
  cusp forms}, Duke Math. J. \textbf{162} (2013), no.~14, 2609--2644.

\bibitem{blomil}
V.~Blomer and D.~Mili\'cevi\'c, \emph{The second moment of twisted modular
  {$L$}-functions}, Geom. Funct. Anal. \textbf{25} (2015), no.~2, 453--516.

\bibitem{butkha}
J.~Buttcane and R.~Khan, \emph{On the fourth moment of {H}ecke-{M}aass forms
  and the random wave conjecture}, Compos. Math. \textbf{153} (2017), no.~7,
  1479--1511.

\bibitem{desiwa}
J.-M. Deshouillers and H.~Iwaniec, \emph{Kloosterman sums and {F}ourier
  coefficients of cusp forms}, Invent. Math. \textbf{70} (1982), no.~2,
  219--288.

\bibitem{golhoflie}
J.~Hoffstein and P.~Lockhart, \emph{Coefficients of {M}aass forms and the
  {S}iegel zero}, Ann. of Math. (2) \textbf{140} (1994), no.~1, 161--181, With
  an appendix by D. Goldfeld, J. Hoffstein and D. Lieman.

\bibitem{ivi}
A.~Ivi\'c, \emph{On sums of {H}ecke series in short intervals}, J. Th\'eor.
  Nombres Bordeaux \textbf{13} (2001), no.~2, 453--468.

\bibitem{iwakow}
H.~Iwaniec and E.~Kowalski, \emph{Analytic number theory}, American
  Mathematical Society Colloquium Publications, vol.~53, American Mathematical
  Society, Providence, RI, 2004.

\bibitem{ils}
H.~Iwaniec, W.~Luo, and P.~Sarnak, \emph{Low lying zeros of families of
  {$L$}-functions}, Inst. Hautes \'Etudes Sci. Publ. Math. (2000), no.~91,
  55--131 (2001).

\bibitem{jut2}
M.~Jutila, \emph{The fourth moment of central values of {H}ecke series}, Number
  theory ({T}urku, 1999), de Gruyter, Berlin, 2001, pp.~167--177.

\bibitem{jut}
\bysame, \emph{The twelfth moment of central values of {H}ecke series}, J.
  Number Theory \textbf{108} (2004), no.~1, 157--168.

\bibitem{kim}
H.~H. Kim, \emph{Functoriality for the exterior square of {${\rm GL}_4$} and
  the symmetric fourth of {${\rm GL}_2$}}, J. Amer. Math. Soc. \textbf{16}
  (2003), no.~1, 139--183, With appendix 1 by D. Ramakrishnan and appendix 2 by
  H. Kim and P. Sarnak.

\bibitem{kiryou}
E.~M. Kiral and M.~P. Young, \emph{The fifth moment of modular $l$-functions},
  preprint, arXiv:1701.07507.

\bibitem{kohzag}
W.~Kohnen and D.~Zagier, \emph{Values of {$L$}-series of modular forms at the
  center of the critical strip}, Invent. Math. \textbf{64} (1981), no.~2,
  175--198.

\bibitem{li}
X.~Li, \emph{Bounds for {${\rm GL}(3)\times {\rm GL}(2)$} {$L$}-functions and
  {${\rm GL}(3)$} {$L$}-functions}, Ann. of Math. (2) \textbf{173} (2011),
  no.~1, 301--336.

\bibitem{pen}
Z.~Peng, \emph{Zeros and central values of automorphic {L}-functions}, ProQuest
  LLC, Ann Arbor, MI, 2001, Thesis (Ph.D.)--Princeton University.

\end{thebibliography}

\end{document}